\documentclass[MSNbibl,number,citesort,seceqn,dvips]{arxbj}
\usepackage{upgreek}
\usepackage{graphicx}
\aid{0}
\volume{20}
\issue{3}
\pubyear{2014}
\firstpage{1097}
\lastpage{1125}
\doi{10.3150/13-BEJ516} 

\makeatletter
\newcommand{\RMo}{\mathrm{o}}

\newcommand{\xrightarrow}[1]{\stackrel{#1}{\longrightarrow}}

\newtheorem{theorem}{Theorem}[section]
\newtheorem{lemma}[theorem]{Lemma}
\newtheorem{proposition}[theorem]{Proposition}
\newtheorem{corollary}[theorem]{Corollary}

\newcommand{\Var}{\operatorname{Var}}
\newcommand{\Cov}{\operatorname{Cov}}

\newcommand{\E}{\mathrm{E}}
\newcommand{\ZZ}{\mathbb{Z}}
\newcommand{\RR}{\mathbb{R}}
\newcommand{\NN}{\mathbb{N}}
\newcommand{\bz}{z} 
\newcommand{\bh}{h} 

\newcommand{\bX}{\mathbf{X}}
\newcommand{\bY}{\mathbf{Y}}
\newcommand{\bC}{\mathbf{C}}

\renewcommand{\div}{\operatorname{div}}
\newcommand{\dd}{\,\mathrm d}

\newcommand{\Par}{{\theta}}
\newcommand{\ParE}{{\widehat{\theta}_n}}

\makeatother

\begin{document}
\begin{frontmatter}

\title{Variational approach for spatial point process intensity estimation}
\runtitle{Spatial point processensity estimation}

\begin{aug}
\author[1]{\inits{J.-F.}\fnms{Jean-Fran\c{c}ois} \snm{Coeurjolly}\corref{}\thanksref{1}\ead[label=e1]{Jean-Francois.Coeurjolly@upmf-grenoble.fr}} \and
\author[2]{\inits{J.}\fnms{Jesper} \snm{M\O{}ller}\thanksref{2}\ead[label=e2]{jm@math.aau.dk}}
\runauthor{J.-F. Coeurjolly and J. M\o{}ller} 
\address[1]{Laboratory Jean Kuntzmann,
Grenoble Alpes University,
51 Rue des Math\'{e}matiques,
BP 53 - 38041 Grenoble Cedex 9,
France. \printead{e1}}
\address[2]{Department of Mathematical Sciences,
Aalborg University,
Fredrik Bajers Vej 7G,
DK-9220 Aalborg,
Denmark. \printead{e2}}
\end{aug}

\received{\smonth{9} \syear{2012}}
\revised{\smonth{1} \syear{2013}}

%
\begin{abstract}
We introduce a new variational estimator for the intensity function of
an inhomogeneous spatial point process with points in the
$d$-dimensional Euclidean space and observed within a bounded region.
The variational estimator applies in a simple and general setting when
the intensity function is assumed to be of log-linear form
$\beta+\Par^\top z(u)$ where $z$ is a spatial covariate function and
the focus is on estimating $\Par$. The variational estimator is very
simple to implement and quicker than alternative estimation procedures.
We establish its strong consistency and asymptotic normality. We also
discuss its finite-sample properties in comparison with the maximum
first order composite likelihood estimator when considering various
inhomogeneous spatial point process models and dimensions as well as
settings were $z$ is completely or only partially known.
\end{abstract}

%

\end{frontmatter}

\section{Introduction}\label{secintro}

Intensity estimation for spatial point processes is of fundamental
importance in many applications,
see, for example, Diggle \cite{diggle03}, M{\o}ller and
Waagepetersen \cite{moellerwaagepetersen07},
Illian \textit{et al.} \cite{illianetal08}, Baddeley~\cite
{baddeley10}, and Diggle \cite{diggle10}. While
maximum likelihood and Bayesian methods are feasible for parametric
Poisson point process models
(Berman and Turner \cite{bermanturner92}), computationally intensive Markov
chain Monte Carlo methods are needed otherwise
(M{\o}ller and Waagepetersen \cite{moellerwaagepetersen04}).
The Poisson likelihood has been used for intensity estimation in
non-Poisson models (Schoenberg \cite{schoenberg05}, Guan and Shen
\cite{guanshen10}) where it
can be viewed as a composite likelihood based on the intensity function
(M{\o}ller and Waagepetersen \cite{moellerwaagepetersen07} and
Waagepetersen \cite{waagepetersen07});
we refer to this as a ``first
order composite likelihood''.
For Cox and Poisson cluster point processes, which form
major classes of point process models for clustering or
aggregation (Stoyan, Kendall and Mecke~\cite
{stoyankendallmecke95}), the first and second
order moment properties as expressed
by
the intensity function $\rho$ and pair correlation function $g$ are
often of
an explicit form, and this has led to the development of
estimation procedures based on combinations of first and second order
composite likelihoods and minimum contrast estimation procedures
(Guan \cite{guan06}, M{\o}ller and Waagepetersen \cite
{moellerwaagepetersen07},
Waagepetersen \cite{waagepetersen07}) and to refinements of such methods
(Guan and Shen \cite{guanshen10}, Guan, Jalilian and Waagepetersen
\cite{waagepetersenguanjalilian11}). For
Gibbs point processes, which form
a major class of point process models for repulsiveness,
the
(Papangelou) conditional intensity is of explicit form and has been
used for developing maximum\vadjust{\goodbreak} pseudo-likelihood estimators (Besag \cite
{besag77},
Jensen and M{\o}ller \cite{jensenmoeller91}, Baddeley and Turner
\cite{baddeleyturner00}) and variational
estimators (Baddeley and Dereudre \cite{baddeleydereudre12}).
However, in general for Gibbs
point processes,
the moment
properties are not expressible in closed form and it is therefore hard to
estimate the intensity function.

The present paper considers a new variational
estimator for the intensity function of a
spatial point process $\bX$, with points in the $d$-dimensional
Euclidean space $\RR^d$ and observed within a bounded region
$W\subset\RR^d$.
It is to some extent derived along similar
lines as the variational
estimator based on the conditional intensity
(Baddeley and Dereudre~\cite{baddeleydereudre12}),
which in turn is a counterpart
of the variational estimator for Markov random fields
(Almeida and Gidas \cite{almeidagidas93}). However, our variational
estimator applies in a much simpler and general setting. 
In analogy with the exponential form of the
conditional intensity considered in Baddeley and Dereudre \cite
{baddeleydereudre12}, we
assume that $\bX$ has a log-linear intensity
function
%
\begin{equation}
\label{erho} \rho(u)=\exp\bigl(\beta+\Par^\top\bz(u) \bigr),\qquad
u\in\RR^d.
\end{equation}
Here $\beta$ is a real parameter, $\Par$ is a real $p$-dimensional
parameter and $\Par^\top$ is its transpose, $\bz$ is a real
$p$-dimensional function defined on
$\RR^d$ and referred to as the covariate function, and we view $\Par$
and $\bz(u)$ as column vectors.
A log-linear intensity
function is often assumed for Poisson point processes (where it is the
canonical link) and for Cox processes (see
M{\o}ller and Waagepetersen \cite{moellerwaagepetersen07} and the
references therein), while
for Gibbs point process models
it is hard to exhibit a model with intensity function of the
log-linear form.
Further details are given in Sections \ref{secset}--\ref{seccde}.

As the variational estimator in Baddeley and Dereudre \cite
{baddeleydereudre12}, our
variational estimator concerns~$\theta$, while $\beta$ is treated as a
nuisance parameter which is not estimated. Our
variational estimator is simple to implement,
it requires only the computation of
the solution of a system of $p$ linear
equations involving
certain sums over
the points of $\bX$ falling in $W$, and it is
quicker to use than the other estimation methods mentioned
above.
Moreover, our
variational estimator is expressible in closed form while
the maximum likelihood estimator for the Poisson
likelihood and the maximum first order composite likelihood
estimator for non-Poisson models are not expressible in
closed form and the profile
likelihood for $\theta$ involves the computation (or approximation) of
$d(1+p/2)(p+1)$ integrals. On the one hand,
as for the approach based on first
order composite likelihoods, an advantage of our
variational estimator is its flexibility, since
apart from (\ref{erho}) and a few mild assumptions on $\bz$,
we do not make any further assumptions. In particular, we do not
require that
$\bX$ is a grand canonical Gibbs process as assumed
in Baddeley and Dereudre \cite{baddeleydereudre12}. On the other hand,
a possible disadvantage of our
variational approach is a loss in efficiency,
since we do not take into account spatial
correlation, for example, through the modelling of
the pair correlation function as in Guan and Shen \cite{guanshen10} and
Guan, Jalilian and Waagepetersen \cite
{waagepetersenguanjalilian11}, or interaction, for example, through
the modelling of the conditional intensity function as in
Baddeley and Dereudre \cite{baddeleydereudre12}.


The paper is organized as follows. Section \ref{secset}
presents our general setting.
Section \ref{seccde} specifies our variational estimator,
establishes its asymptotic properties, and discusses the conditions we
impose.
Section \ref{secsim} reports on a simulation study of
the finite-sample properties of our
variational estimator and the maximum first order composite likelihood
estimator for various inhomogeneous spatial point process models
in the planar case $d=2$ as well as higher dimensions and when $z$ is
known on an observation window as well as when $z$ is
known only on a finite set of locations.
The technical proofs of our results
are deferred to Appendix \ref{appA}. Finally,
Appendix \ref{appmodel2} illustrates the simplicity of
our variational estimator and the flexibility of the conditions
given in Section \ref{seccde}.



\section{Preliminaries}\label{secset}

This section introduces the assumptions and notation used throughout
this paper.

Let $W\subset\RR^d$ be a compact set
of positive Lebesgue measure $|W|$. It will play the role of an
observation window. Without any danger of confusion,
we also use the notation $|A|$
for the cardinality of a countable set $A$, and
$|u|=\max\{|u_i|\dvt i=1,\ldots,d\}$ for the maximum
norm of a point $u=(u_1,\ldots,u_d)\in\RR^d$. Further, we let
$\|u\|$ denote the Euclidean norm for a point $u\in\RR^d$, and $\|A\|
={\sup_{\|u\|=1}}|Au|$
the supremum norm for a square matrix $A$, that is, its numerically
largest (right) eigenvalue. Moreover,
for any real $p$-dimensional function $k$ defined on $\RR^d$, we let
%
\begin{equation}
\label{enorm} \|k\|_\infty= \sup_{u \in\RR^d} \bigl\| k(u) \bigr\|.
\end{equation}

Let $\bX$ be a spatial point process on $\RR^d$, which
we view as a random
locally finite subset of~$\RR^d$. Let $\bX_W=\bX\cap W$.
Then the number of points in
$\bX_W$ is finite; we denote this number by
$N(W)=n(\bX_W)=|\bX_W|$; and
a realization of $\bX_W$ is of the form
$\mathbf{x}=\{x_1,\ldots,x_n\}\subset W$, where
$n=n(\mathbf{x})$ and
$0\le n<\infty$. If $n=0$, then $\mathbf{x}=\varnothing$ is the empty point
pattern in $W$. For further background material and measure theoretical details
on spatial point process, see, for example, Daley and Vere-Jones \cite
{daleyvere-jones03} and
M{\o}ller and Waagepetersen \cite{moellerwaagepetersen04}.

We assume that $\bX$ has a locally integrable intensity function
$\rho$. By Campbell's theorem (see, e.g., M{\o}ller and
Waagepetersen \cite{moellerwaagepetersen04}),
for any real Borel function $k$ defined on $\RR^d$ such that $k\rho$
is absolutely
integrable (with respect to the Lebesgue measure on $\RR^d$),
%
\begin{equation}
\label{ecampbell} \E\sum_{u\in\bX}k(u)=\int k(u)\rho(u)\,
\mathrm du.
\end{equation}
Furthermore,
for any
integer $n\ge1$, $\bX$ is said to have an $n$th order product density
$\rho^{(n)}$ if this is a non-negative Borel function on $\RR^{dn}$
such that for all non-negative Borel functions $k$ defined on~$\RR^{dn}$,
%
\begin{equation}
\label{ecampbell2} \E\sum^{\neq}_{u_1,\ldots,u_n\in\bX}k(u_1,\ldots,u_n)=\int\cdots\int k(u_1,\ldots,u_n)
\rho^{(n)}(u_1,\ldots,u_n) \,\mathrm
du_1\cdots\mathrm du_n,
\end{equation}
where
the $\neq$ over the summation sign means that $u_1,\ldots,u_n$ are
pairwise distinct. Note that $\rho=\rho^{(1)}$.

Throughout this paper except in
Section \ref{secbasic}, we assume that $\rho$
is of the log-linear form
(\ref{erho}),
where we view $\Par$
and $\bz(u)$
as $p$-dimensional
column vectors.

As for vectors, transposition of a matrix $A$ is denoted
$A^\top$. For convenience, we, for example, write $(\beta,\Par)$ when
we more
precisely mean the $(p+1)$-dimensional
column vector $(\beta,\Par^\top)^\top$. If $A$ is a square matrix, we
write $A\ge0$ if $A$ is positive
semi-definite, and $A>0$ if $A$ is (strictly) positive definite. When
$A$ and $B$ are square matrices of the same size, we write
$A\ge B$ if \mbox{$A-B\ge0$}. 

For $k=0,1,\ldots\,$, denote $\mathcal{C}_{d,p}^k$ the class of
$k$-times continuous differentiable
real $p$-dimensional
functions defined on $\RR^d$.
For $h\in\mathcal{C}_{d,1}^1$, denote its gradient
\[
\nabla h(u)= \biggl(\frac{\partial h}{\partial
u_1}(u),\ldots,\frac{\partial h}{\partial u_d}(u)
\biggr)^\top,\qquad u=(u_1,\ldots,u_d)^\top
\in\RR^d
\]
and define the divergence operator $\div$ on $\mathcal{C}_{d,1}^1$ by
\[
\div h (u)= \frac{\partial h}{\partial u_1}(u) + \cdots+ \frac
{\partial h}{\partial u_d}(u),\qquad
u=(u_1,\ldots,u_d)^\top\in\RR^d.
\]
Furthermore, for $h=(h_1,\ldots,h_p)^\top\in\mathcal{C}_{d,p}^1$,
define the divergence operator $\div$ on $\mathcal{C}_{d,p}^1$ by
\[
\div\bh(u)= \bigl(\div h_1(u),\ldots, \div h_p(u)
\bigr)^\top,\qquad u\in\RR^d.
\]
%
If $\bz\in\mathcal{C}_{d,p}^1$, then by (\ref{erho})
%
\begin{equation}
\label{ediv-log-rho} \div\log\rho(u)=\Par^\top\div\bz(u)= \div
\bz(u)^\top\Par,\qquad u\in\RR^d.
\end{equation}

Finally, we recall the classical definition of mixing coefficients
(see, e.g., Politis, Paparoditis and Romano \cite{politis98}): for
$j,k \in\NN\cup\{\infty\}$ and
$m\geq1$, define
\begin{eqnarray*}
\alpha_{j,k}(m)&=&\sup\bigl\{ \bigl|P(A\cap B) - P(A)P(B)\bigr|\dvt  A\in\mathcal{F}(
\Lambda_1), B\in\mathcal{F}(\Lambda_2),
\\
&&\hspace*{20.6pt} \Lambda_1\in\mathcal B\bigl(\RR^d\bigr),
\Lambda_2 \in\mathcal B\bigl(\RR^d\bigr), |
\Lambda_1|\leq j, |\Lambda_2|\leq k, d(
\Lambda_1,\Lambda_2)\geq m \bigr\},
\end{eqnarray*}
where $\mathcal{F}(\Lambda_i)$ is the $\sigma$-algebra generated by
$X\cap\Lambda_i$, $i=1,2$, $d(\Lambda_1,\Lambda_2)$ is the
minimal distance between the sets $\Lambda_1$ and $\Lambda_2$, and
$\mathcal B(\RR^d)$ denotes the class of Borel sets in $\RR^d$.


\section{The variational estimator}\label{seccde}


Section \ref{secbasic}
establishes an identity which together with (\ref{ediv-log-rho})
is used in Section \ref{secnew-est} for deriving an
unbiased estimating equation which only involves $\Par$, the
parameter of interest, and from which our variational estimator is
derived. Section \ref{secasym-cde} discusses the asymptotic
properties of the variational estimator.

\subsection{Basic identities}\label{secbasic}

This section establishes some basic identities for a spatial point
process $\bX$ defined on $\RR^d$ and having a locally integrable
intensity function
$\rho$ which is not necessarily of the log-linear form (\ref{erho}).
The results will be used later when defining our variational estimator.

Consider a real Borel function $h$ defined
on $\RR^d$ and let $f(u)=\rho(u)|h(u)|$.
For $n=1,2,\ldots\,$, let 
$E_n^d=[-n,n]^d$ and
\[
\mu_{n}(f)=\max\bigl\{\mu_{n,j}(f)\dvt j=1,\ldots,d\bigr\}
\]
with
\begin{eqnarray*}
\mu_{n,j}(f)&=&\int_{E_n^{d-1}} f(u_1,\ldots,u_{j-1},-n,u_{j+1},\ldots,u_d)\,\mathrm
du_1\cdots\mathrm du_{j-1}\,\mathrm du_{j+1}\cdots
\mathrm du_n
\\
&&{}+\int_{E_n^{d-1}} f(u_1,\ldots,u_{j-1},n,u_{j+1},\ldots,u_d)\,\mathrm du_1\cdots\mathrm du_{j-1}\,
\mathrm du_{j+1}\cdots\mathrm du_n
\end{eqnarray*}
provided the integrals exist.
Note that $\mu_{n}(f)$ depends only on
the behaviour of $f$ on 
the boundary of
$E_n^d$.

\begin{proposition}\label{prop4} Suppose that $h,\rho\in\mathcal
C^1_{d,1}$ such that $\lim_{n\rightarrow\infty}\mu_n(\rho|h|)=0$ and
for $j=1,\ldots,d$, the function $h(u)\,\partial\rho(u)/\partial u_j$ is
absolutely
integrable. Then the following relations hold where the
mean values exist and are
finite:
%
\begin{equation}
\label{ei1} \E\sum_{u\in\bX} h(u) \nabla\log\bigl(\rho(u)
\bigr) = - \E\sum_{u \in\bX} \nabla h(u)
\end{equation}
and
%
\begin{equation}
\label{ei2} \E\sum_{u \in\bX} h(u) \div\log\bigl(\rho(u)
\bigr) =- \E\sum_{u \in
\bX} \div h(u).
\end{equation}
\end{proposition}

\begin{pf}
For
$j=1,\ldots,d$ and $u=(u_1,\ldots,u_d)^\top\in\RR^d$, Campbell's theorem
(\ref{ecampbell}) and the assumption that
$h(u)\,\partial\rho(u)/\partial u_j$ is absolutely
integrable
imply that
\[
\E\biggl( \sum_{u \in\bX} h(u) \nabla\log\bigl(\rho(u)
\bigr) \biggr)_j = \int h(u) \,\frac{\partial\rho}{\partial u_j}(u)\,
\mathrm du
\]
exist. Thereby,
\begin{eqnarray*}
&&
\E\biggl( \sum_{u \in\bX} h(u) \nabla\log\bigl(\rho(u)
\bigr) \biggr)_j \\
&&\quad=\lim_{n\to\infty} \int
_{E_n^d} h(u)\,\frac{\partial\rho}{\partial u_j} (u)\,\mathrm du
\\
&&\quad= \lim_{n\to\infty} \int_{E_n^{d-1}} \biggl( \bigl[
\rho(u)h(u) \bigr]_{u_j=-n}^{u_j=n} - \int_{-n}^n
\rho(u) \,\frac{\partial h}{\partial u_j} (u)\,\mathrm du_j \biggr)\,
\mathrm
du_1 \cdots\mathrm du_{j-1}\,\mathrm du_{j+1}
\cdots\mathrm du_n
\\
&&\quad= -\lim_{n\to\infty} \int_{E_n^d} \frac{\partial h}{\partial u_j}
(u) \rho(u)\,\mathrm du,
\end{eqnarray*}
where the first identity follows 
from the dominated convergence
theorem, the second from Fubini's theorem and
integration by parts, and the third from Fubini's theorem and the
assumption that $\lim_{n\rightarrow\infty}\mu_n(\rho|h|)=0$, since
\[
\biggl|\int_{E_n^{d-1}} \bigl[ \rho(u)h(u) \bigr]_{u_j=-n}^{u_j=n}
\biggr| \le\mu_{n,j}\bigl(\rho|h|\bigr)\le\mu_n\bigl(\rho|h|\bigr).
\]
Hence, using first the dominated convergence
theorem and second Campbell's theorem,
\[
\E\biggl( \sum_{u \in\bX} h(u) \nabla\log\bigl(\rho(u)
\bigr) \biggr)_j = 
-\int\frac{\partial h}{\partial u_j} (u) \rho(u)\,
\mathrm du
=-\E\biggl( \sum_{u \in\bX}
\nabla h(u) \biggr)_j
\]
whereby (\ref{ei1}) is verified and the mean values in (\ref{ei1})
are seen to exist and are finite.
Finally, (\ref{ei1}) implies (\ref{ei2}) where the mean values
exist and are finite.
\end{pf}

Proposition \ref{prop4}
becomes useful when $\rho$ is of the log-linear form (\ref{erho}):
if we omit the
expectation signs in (\ref{ei1})--(\ref{ei3}), we obtain unbiased
estimating equations, where
(\ref{ei1}) gives
a linear system of $p$ vectorial equation in dimension $d$,
while (\ref{ei3}) gives a linear system of $p$ one-dimensional
equations for the estimation of the $p$-dimensional parameter
$\theta$;
the latter system is simply obtained by summing over
the $d$ equations in each vectorial equation.
A similar reduction of equations is obtained in Baddeley and Dereudre
\cite{baddeleydereudre12}.

The conditions and the last result in
Proposition \ref{prop4}
simplify
as follows when $h$ vanishes outside $W$.

\begin{corollary}\label{cor1} Suppose that $h,\rho\in\mathcal
C^1_{d,1}$ such that $h(u)=0$ whenever $u\notin
W$. Then
%
\begin{equation}
\label{ei2a} \E\sum_{u \in\bX_W} h(u) \div\log\bigl(\rho(u)
\bigr) =- \E\sum_{u \in
\bX
_W} \div h(u).
\end{equation}
\end{corollary}

\subsection{The variational estimator}\label{secnew-est}

Henceforth we consider the case of the log-linear intensity function
(\ref{erho}), assuming that the parameter space for
$(\beta,\Par)$ is $\RR\times\RR^p$.
We specify below our variational estimator in terms of a
$p$-dimensional
real test
function
\[
\bh=(h_1,\ldots,h_p)^\top
\]
defined on $\RR^d$. The test function
is required not to depend on $(\beta,\Par)$ and to satisfy certain
smoothness conditions. The specific choice of test functions is
discussed at the end of Section \ref{secnewsec}.

In the present section, to stress that the expectation of a
functional $f$ of $\bX$ depends on $(\beta,\Par)$, we write this as
$\E_{\beta,\Par}f(\bX)$. Furthermore,
define the $p\times p$ matrix
\[
A(\bX_W)= \sum_{u \in\bX_{W}} h(u) \div
z(u)^\top
\]
and the $p$-dimensional column vector
\[
b(\bX_W)=\sum_{u \in\bX_{W}} \div h(u).
\]

\subsubsection{Estimating equation and definition of the variational
estimator}\label{secesteqdef}

We consider first the case where the test function $h$ vanishes outside $W$.

\begin{corollary} \label{corest}
Suppose that $\bh,\bz\in\mathcal
C^1_{d,p}$ such that
%
\begin{equation}
\label{einappr}
\mbox{$\bh(u)=0$}\qquad \mbox{whenever $u\notin W$.}
\end{equation}
Then, for any $(\beta,\Par)\in\RR\times\RR^p$,
%
\begin{equation}
\label{ei3} \E_{\beta,\Par} A(\bX_W) \Par= - \E_{\beta,\Par} b(
\bX_W) .
\end{equation}
\end{corollary}

\begin{pf} The conditions of Corollary \ref{cor1} are easily seen
to be satisfied. Hence combining (\ref{ediv-log-rho}) and
(\ref{ei2a}) we obtain (\ref{ei3}).
\end{pf}

Several remarks are in order.

Note that (\ref{ei3}) is a linear system of $p$
equations for the $p$-dimensional parameter $\theta$.
Under the conditions in Corollary \ref{corest},
(\ref{ei3})
leads to the unbiased estimating equation
%
\begin{equation}
\label{eunbiasedesteq} A(\bX_W) \theta=-b(\bX_W).
\end{equation}
Theorem \ref{thmconv} below establishes that under certain conditions,
where we
do not necessarily require $h$ to vanish outside $W$,
(\ref{eunbiasedesteq})
is an asymptotically
unbiased estimating equation as $W$ extends to~$\RR^d$.

In the sequel we therefore
do not necessarily
assume (\ref{einappr}).
For instance, when $\div z(u)$ does not vanish outside $W$,
we may consider either $h(u)= \div z(u)$ or $h(u)= \eta_W(u) \div
z(u)$, where $\eta_W$ is a smooth function which vanishes outside $W$.
In the latter case, (\ref{eunbiasedesteq}) is an unbiased
estimating equation, while in the former case it is
an asymptotically
unbiased estimating equation (under the conditions imposed in
Theorem \ref{thmconv}).

When (\ref{eunbiasedesteq})
is an (asymptotically)
unbiased estimating equation and $A(\bX_W)$ is invertible,
we \textit{define the variational
estimator} by
%
\begin{equation}
\label{edaggert} \widehat\theta=-A(\bX_W)^{-1}b(
\bX_W).
\end{equation}
Theorem \ref{thmconv} below establishes under certain conditions
the invertibility of $A(\bX_W)$ and the strong consistency and
asymptotic normality of $\widehat\theta$ as $W$ extends to~$\RR^d$.

Finally, if $h$ is allowed to depend on $\theta$,
(\ref{eunbiasedesteq}) still provides an unbiased estimating equation
but the closed form expression (\ref{edaggert}) only applies when $h$
is not depending on $\theta$ (as assumed in this paper).

\subsubsection{Choice of test function}\label{secnewsec}

The choice of test function should take into consideration the
conditions introduced later in Section \ref{secconditions}.
The test functions below are defined in terms of the covariate
function so that it is possible to check these conditions as discussed
in Section \ref{secdisccond}.

Interesting choices of the test function include:
\begin{itemize}
\item$h(u)= \div z(u)$ and the corresponding modification $h(u)=\eta
_W(u)\div z(u)$,
\item$h(u)=z(u)$ and the corresponding modification $h(u)=\eta_W(u)z(u)$.
\end{itemize}
In the first case, $A(\bX_W)$ becomes a covariance matrix. For example,
if $h(u)= \div z(u)$, then
\[
A(\bX_W)= \sum_{u \in\bX_{W}} \div z(u) \div
z(u)^\top
\]
%
is invertible if and only if
$A(\bX_W)>0$, meaning that if $\bX_W=\{x_1,\ldots,x_n\}$ is observed,
then the
$p\times n$ matrix with columns
$\div\bz(x_1),\ldots,\div\bz(x_n)$ has rank $p$. In the latter case,
$A(\bX_W)$ is in general not symmetric and we avoid the
calculation of $\div\div z(u)$.

\subsubsection{Choice of smoothing function}

\begin{figure}

\includegraphics[scale=0.98]{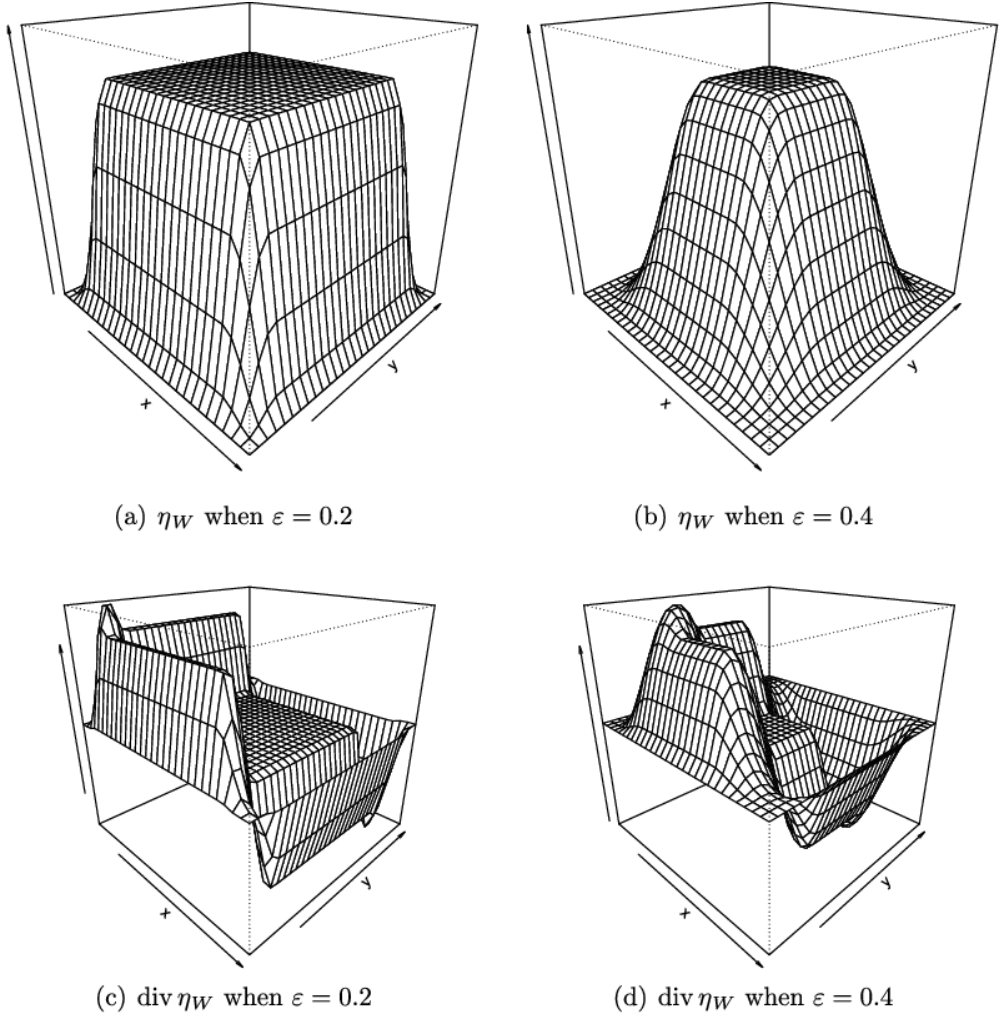}

\caption{Plots of the functions
$\eta_W=\chi_W*\varphi_\varepsilon$ and $\div\eta_W$ when
$W=[-1,1]^2$ and $\varepsilon=0.2,0.4$.}
\label{figmoll}\vspace*{-3pt}
\end{figure}

We let henceforth the
smoothing function
$\eta_{W}$ depend
on a user-specified parameter $\varepsilon>0$ and define it
as the convolution
%
\begin{equation}
\label{eeta} \eta_W(u)=\chi_{W_{\ominus\varepsilon}}*\varphi_\varepsilon
(u)=\int\mathbf{1}(u-v\in W_{\ominus\varepsilon})\varphi_\varepsilon
(v)\dd v,\qquad u\in\RR^d,
\end{equation}
where the notation means the following:
\[
{W}_{\ominus\varepsilon}=\bigl\{u\in W\dvt b(u,\varepsilon)\subseteq W\bigr
\}
\]
is the observation window eroded by the $d$-dimensional closed ball
$b(u,\varepsilon)$ centered at $u$ and with radius $\varepsilon$;
$\chi_{W_{\ominus\varepsilon}}(\cdot)=\mathbf{1}(\cdot\in W_{\ominus
\varepsilon})$ is the indicator function on $W_{\ominus\varepsilon}$;
and
\[
\varphi_\varepsilon(u)=\varepsilon^{-d}\varphi(u/\varepsilon),\qquad
u\in\RR^d,
\]
where
\[
\varphi(u) = c \exp\biggl({-\frac{1}{1-\|u\|^2}} \biggr) \mathbf{1}\bigl(
\|u\|\leq1\bigr),\qquad u\in\RR^d,
\]
where
$c$ is a normalizing constant such that $\varphi$ is a density
function ($c \approx2.143$ when $d=2$).
Figure \ref{figmoll} shows the function $\eta_W$
and its divergence when $W=[-1,1]^2$, $\varepsilon=0.2$, and
$\varepsilon=0.4$.
The
construction (\ref{eeta})
is quite standard in distribution theory when functions are
regularized and it can be found, though in a slightly
different form, in H{\"o}rmander (\cite{hormander03}, Theorem 1.4.1,
page~25).

It is easily checked that $\varphi_\varepsilon\in\mathcal
C^\infty_{d,1}$, and so $\eta_W \in\mathcal
C^\infty_{d,1}$.
Note that
%
\begin{equation}
\label{eeta2} 0\le\eta_W\le1,\qquad \eta_W(u)=1\qquad\mbox{if $u
\in W_{\ominus2\varepsilon}$},\qquad \eta_W(u)=0\qquad\mbox{if $u\notin W$}.
\end{equation}

The following lemma states some properties for test
functions of the modified form
$h(u)= \eta_W(u) k(u)$, where we let
$\kappa= \int_{\mathcal B(0,1)} |{\div\varphi(v)}|\dd v$; if $d=2$ then
$\kappa\approx1.256 $.


\begin{lemma}\label{lemhepsilon}
Let\vspace*{1pt} $k\in\mathcal{C}^1_{d,p}$ and $h(u)= \eta_W(u) k(u)$ where
$\eta_W$ is given by (\ref{eeta}). Then $h\in\mathcal{C}^1_{d,p}$ and
its support is included in $W$. Further,\vadjust{\goodbreak}
$h$ respective $\div h$ agrees with $k$ respective $\div k$ on
${W}_{\ominus2\varepsilon}$. Moreover, for any $u \in{W}$,
%
\begin{equation}
\label{ehk2} \bigl\| h(u) \bigr\| \le\bigl\| k(u) \bigr\|,\qquad
\bigl\| \div h(u) -\div k(u)\bigr\| \leq\bigl\|
\div k(u)
\bigr\|+ \bigl\| k(u) \bigr\| { \kappa}/{\varepsilon}.
\end{equation}
\end{lemma}

\begin{pf}
We have $h\in\mathcal{C}^1_{d,p}$ since $k\in\mathcal{C}^1_{d,p}$ and
$\eta_W \in\mathcal C^\infty_{d,1}$, and the
support of $h$ is included in $W$ since $\eta_W(u)=0$ if $u\notin W$.
From the last two statements of (\ref{eeta2}), we obtain
that $\div h(u)$ agrees with $\div k(u)$ on $W_{\ominus2\varepsilon}$.
The first inequality in (\ref{ehk2}) follows immediately from the definition
of~$h$, since $\| h(u) \| = \| \eta_W(u) k(u) \|\le\| k(u) \|$.
Recall that $(f*g)'=f*g'$ if $g\in C^1_{d,p}$ has compact
support and $f$ is Lebesgue integrable on $\RR^d$, where in our case
we let $f=\chi_{W_{\ominus\varepsilon}}$ and $g=\varphi_\varepsilon$.
Therefore and
since $\div\varphi_\varepsilon=(\div\varphi)/\varepsilon\in
\mathcal{C}^\infty_{d,1}$, for any $u\in W$, we have
\begin{eqnarray*}
\div h(u) &=& \eta_W(u)\div k(u) + k(u) (
\chi_{{W}_{\ominus
\varepsilon}} * \div\varphi_\varepsilon) (u)
\\
&=& \eta_W(u) \div k(u) + \frac1\varepsilon k(u) ( \chi
_{W_{\ominus
\varepsilon}} * \div\varphi) (u).
\end{eqnarray*}
Thereby, the second inequality in (\ref{ehk2}) follows
from a straightforward calculation using again the fact that $\eta
_W(u)\le1$.
\end{pf}

\subsection{Asymptotic results}\label{secasym-cde}\vspace*{-2pt}

In this section, we present asymptotic results for the variational
estimator when considering a sequence of observation windows $W=W_n$,
$n=1,2,\ldots\,$, which expands to $\RR^d$ as
$n\to\infty$, and a corresponding sequence of test functions $h=h^{(n)}$,
$n=1,2,\ldots\,$. Corresponding to the two cases of test functions
considered in
Section \ref{secesteqdef}, we consider the following two cases:
\begin{enumerate} 
\item[(A)] either $h^{(n)}=k$ does not depend on $n$, 
\item[(B)] or $h^{(n)}(u)= \eta_{W_n}(u)k(u)$, where
$\eta_{W_n}$ is given by (\ref{eeta}).\vadjust{\goodbreak}
\end{enumerate}

\subsubsection{Conditions}\label{secconditions}

Our asymptotic results require the following conditions.

We restrict attention to the spatial case $d\ge2$ (this is mainly for
technical reasons as explained in Section \ref{secmainresult}).
We suppress in the notation that
the intensity $\rho$
and the higher order product densities $\rho^{(2)},\rho^{(3)},\ldots$
depend on the ``true parameters'' $(\beta,\theta)$. Let
%
\begin{equation}
\label{e25} S_n = \int_{W_n} h^{(n)}(u)
\div z(u)^\top\rho(u)\,\mathrm du
\end{equation}
and
%
\begin{equation}
\label{e26} \Sigma_n = \int_{W_n}
f^{(n)}_\theta(u)f^{(n)}_\theta(u)^\top
\rho(u)\dd u + \int_{W_n^2} f^{(n)}_\theta(u_1)f^{(n)}_\theta(u_2)^\top
Q_2(u_1,u_2)\dd u_1\dd
u_2,
\end{equation}
where $Q_2(u_1,u_2)=\rho^{(2)}(u_1,u_2)-\rho(u_1)\rho(u_2)$
(assuming $\rho^{(2)}$ exists) and
\[
f_\theta^{(n)}(u) = h^{(n)}(u) \div
z(u)^\top\theta+ \div h^{(n)}(u),\qquad u\in\RR^d.
\]
It will follow from
the proof of
Theorem \ref{thmconv} below that under the conditions (i)--(vi) stated
below, with probability
one, the integrals in (\ref{e25})--(\ref{e26}) exist and are finite
for all
sufficiently large $n$.

We impose the following conditions, where $o$ denotes the origin of
$\RR^d$:

\begin{enumerate}[(iii)]
\item[(i)] For every $n\geq1$, $W_n=n A=\{na\dvt a\in A\}$,
where $A\subset\RR^d$ is
convex, compact, and contains $o$ in its interior.
\item[(ii)] The test functions $h^{(n)}$, $n=1,2,\ldots\,$, and the
covariate function $z$ are elements of $\mathcal{C}^1_{d,p}$, and
satisfy for some
constant $K>0$,
%
\begin{eqnarray}
\label{estjerne} \| z\|_\infty&\leq& K,\qquad \| {\div z}\|_\infty\leq
K,\nonumber\\[-8pt]\\[-8pt]
\sup_{n\geq1} \bigl\| h^{(n)}\bigr\|_\infty&\leq& K,\qquad \sup
_{n\geq1} \bigl\| \div h^{(n)}\bigr\|_\infty\leq K.\nonumber
\end{eqnarray}
%

\item[(iii)] There exists a $p\times p$ matrix $I_0$ such
that for all sufficiently large $n$, we have
$S_n/|W_n|\geq I_0>0$.

\item[(iv)] There exists an integer $\delta\geq1$ such that for
$k=1,\ldots,2+\delta$, the product density
$\rho^{(k)}$ exists and $\rho^{(k)}\le K'$, where $K'<\infty$ is  a
constant.


\item[(v)] For the strong mixing coefficients (Section \ref{secset}),
we assume that there exists some $\nu> d(2+\delta)/\delta$ such that
$a_{2,\infty}(m)=\mathcal{O}(m^{-\nu})$.

%
\item[(vi)] The second order product density $\rho^{(2)}$ exists, and
there exists a $p\times p$ matrix $I'_0$
such that for all sufficiently large $n$, $\Sigma_n/|W_n|\geq
I'_0>0$.
\end{enumerate}

\subsubsection{Discussion of the conditions}\label{secdisccond}

Some comments on conditions (i)--(vi) are in order.

In general in applications, the observation window has a non-empty
interior. In (i), 
the assumption that $A$ contains $o$ in its interior can be made
without loss of generality; if instead $u$ was an interior point of
$A$, then
(i) could be modified to that any
ball with centre $u$ and radius $r>0$ is contained in $W_n=nA$ for all
sufficiently large $n$. We could also modify (i) to the case where
$|A|>0$ and as $n\rightarrow\infty$ the limit of $W_n=nA$ exists and
is given by $W_\infty$; then in
(\ref{estjerne}) we should redefine
$\|\cdot\|_\infty=\sup_{u\in\RR^d}\|k(u)\|$ (i.e., as defined in
(\ref{enorm})) by $\|\cdot\|_\infty=\sup_{u\in
W_\infty}\|k(u)\|$. For either case, Theorem \ref{thmconv} in
Section \ref{secmainresult}
will remain true, as the proof of the theorem (given in Appendix \ref{appA})
can easily be modified
to cover these cases.

In (ii), for both cases of (A) and (B) and for $k(u)=\div z(u)$,
(\ref{estjerne}) simplifies to
%
\begin{equation}
\label{euseful} \| z\|_\infty\leq K,\qquad \| {\div z}\|_\infty\leq K,\qquad
\|{\div\div z}\|_\infty\le K.
\end{equation}
This follows immediately for the case (A), since then $h^{(n)}=h$ does
not depend on $n$, while in the case (B) where
$h^{(n)}(u)=
\eta_{W_n}(u) k(u)$,
Lemma \ref{lemhepsilon} implies the equivalence of (\ref{estjerne})
and (\ref{euseful}).

Note that in (ii) we
do not require that $h^{(n)}$ vanishes outside $W_n$. Thus, in connection
with the unbiasedness result in
Corollary \ref{corest}, one of the difficulties to prove
Theorem \ref{thmconv} below will be to ``approximate'' $h^{(n)}$ by a
function with support $W_n$, as detailed in Appendix \ref{appA}.

Conditions (iii) and (vi) are spatial average assumptions
like when establishing asymptotic normality of ordinary least square
estimators for linear models.
These conditions
must be checked for each choice of covariate function, since they
depend strongly on $z$. Note that under condition~(ii), for any
$u\in\RR^d$, $\rho(u) \geq\exp(\beta-\|\theta\|_\infty\|z\|
_\infty
)=c>0$. Therefore, condition (iii)
is satisfied if $h^{(n)}(u)\div z(u)^\top\geq0$ for any $u$ and if
$|W_n|^{-1}\int_{W_n} h^{(n)}(u)\div z(u)^\top\dd u\geq I_0$ for all
sufficiently large $n$. In addition, if
$Q_2(u_1,u_2) \geq0$ for any $u_1,u_2\in\RR^d$ (this is discussed above
for \mbox{specific} point process models), then condition (vi) is satisfied if
$|W_n|^{-1}\int_{W_n} f^{(n)}(u)f^{(n)}(u)^\top\dd u \geq I_0^\prime$
for all sufficiently large $n$.

Condition (iv) is not very restrictive. It is fulfilled
for any
Gibbs
point process with a
Papangelou conditional intensity which is uniformly bounded from above
(the so-called local stability
condition, see, e.g., M{\o}ller and Waagepetersen \cite
{moellerwaagepetersen04}), and also for a
log-Gaussian Cox process where the mean and covariance functions of
the underlying Gaussian process are uniformly bounded from above
(see M{\o}ller, Syversveen and Waagepetersen \cite
{moellersyversveenwaagepetersen98} and
M{\o}ller and Waagepetersen \cite{moellerwaagepetersen07}).
Note that the
larger we can choose $\delta$, the weaker becomes condition (v).

Condition (v) combined with (iv) is
also considered in Waagepetersen and Guan \cite
{waagepetersenguan09}, and (iv)--(v) are
inspired by a central limit theorem obtained first
by Bolthausen \cite{bolthausen82} and later
extended to non-stationary random fields in Guyon \cite{guyon91}
and to triangular arrays of non-stationary random fields (which is
the requirement of our setting) in Kar{\'a}csony \cite{karaczony06}.
We underline
that we turned to a
central limit theorem using mixing conditions instead of one using
martingale type assumptions (e.g., Jensen and K{\"u}nsch \cite
{jensenkunsch94})
since for most of models considered in this paper (in particular the
two Cox processes discussed below) the
``martingale'' type assumption is not satisfied. Such an assumption is
more devoted to Gibbs point processes.\vadjust{\goodbreak}

Other papers dealing with asymptotics for estimators based on estimating
equations for spatial point processes (e.g., Guan \cite{guan06},
Guan and Loh \cite{guanloh07}, Guan and Shen \cite{guanshen10},
Guan, Jalilian and Waagepetersen \cite
{waagepetersenguanjalilian11}, Proke{\v{s}}ov{\'a} and Jensen \cite
{prokesovajensen10})
are assuming mixing
properties expressed in terms 
of a different definition of mixing coefficient (see, e.g.,
Equations (5.2)--(5.3) in Proke{\v{s}}ov{\'a} and Jensen \cite
{prokesovajensen10}).
The mixing conditions in these papers are related to
a central limit theorem by Ibragimov and Linnik \cite
{ibramigovlinnik71} obtained
using blocking techniques, and the mixing conditions may
seem slightly less restrictive than our condition (v). However,
rather than our condition (iv), it is assumed in the papers
that the first four
reduced cumulants exist and have finite total variation. In our
opinion, this is an awkward assumption in the case of Gibbs point
processes and many other examples of spatial point
process models, including Cox processes where the first four
cumulants are not (easily) expressible in a closed form (one exception
being log-Gaussian Cox processes).

Condition (v) is also discussed in
(Waagepetersen and Guan \cite{waagepetersenguan09}, Section 3.3 and
Appendix~E) from which we
obtain that (v) is satisfied in, for example, the following cases of a Cox
process~$\bX$.
\begin{itemize}
\item An inhomogeneous
log-Gaussian Cox process (M{\o}ller and Waagepetersen \cite
{moellerwaagepetersen07}):
Let $\bY$ be a Gaussian process
with mean function $m(u)=\beta+\Par^\top\bz(u)-\sigma^2/2$,
$u\in\RR^2$, and a stationary covariance function
$c(u)=\sigma^2r(u)$, $u\in\RR^2$, where
$\sigma^2>0$ is the variance and the
correlation function $r$ decays at a rate faster
than $d+\nu$. This includes the case of the exponential
correlation function which is considered later in
Section \ref{secsimd2}. If
$\bX$ conditional on $\bY$ is a Poisson point process with intensity
function $\exp(\bY)$, then $\bX$ is an inhomogeneous
log-Gaussian Cox process.
\item An inhomogeneous Neyman--Scott process
(M{\o}ller and Waagepetersen \cite{moellerwaagepetersen07}): Let
$\bC$ be a stationary Poisson point process with intensity $\kappa>0$,
and $f_\sigma$ a
density function on $\RR^d$ satisfying
\[
\sup_{w \in[-m/2,m/2]^d} \int_{\RR^d \setminus[-m,m]^d} f_\sigma
(v-w)\dd w =\mathcal O\bigl(m^{-\nu}\bigr).
\]
This includes the case where $f_\sigma$ is the density function of
$\mathcal N(0,\sigma^2I_d)$, that is, the zero-mean isotropic
$d$-dimensional
normal distribution with standard deviation $\sigma>0$;
we consider this case
later in Section \ref{secsimd2}. If
$\bX$ conditional on $\bC$ is a Poisson point process with intensity
function
%
\begin{equation}
\label{einhomSNCP} \exp\bigl(\beta+\Par^\top\bz(u)\bigr)\sum
_{c\in\bC}f_\sigma(u-c)/\kappa,\qquad u\in\RR^2,
\end{equation}
then $\bX$ is an inhomogeneous Neyman--Scott process.
When $f_\sigma$ is the density function of
$\mathcal N(0,\sigma^2I_d)$, we refer to $\bX$ as an inhomogeneous
Thomas process.
\end{itemize}
Note that in any of these cases of Cox processes, $\rho(u)=\exp(\beta
+\Par^\top
\bz(u))$ is indeed an intensity function of the log-linear form
(\ref{erho}) and that for both cases the pair correlation function is greater
than $1$ which implies that $Q_2(u_1,u_2)\geq0$ for any $u_1,u_2 \in
\RR^d$.

Moreover, for Gibbs point processes, (v) may be checked using results
in Heinrich \cite{heinrich92} and Jensen \cite{jensen93}, where in
particular
results for pairwise interaction point processes satisfying a
hard-core type condition may apply. However, as stressed in Section \ref{secintro},
the problem with Gibbs models is
that it is hard to exhibit a model with intensity function of the
log-linear form (\ref{erho}).

Finally, if $\bX$ is a Poisson point process many simplifications
occur. First, for any integer $k\ge1$,
$\rho^{(k)}(u_1,\ldots,u_k)=\rho(u_1)\cdots\rho(u_k)$, and hence (iv)
follows from (ii). Second, since
$\bX_{\Lambda_1}$ and $\bX_{\Lambda_2}$ are
independent whenever $\Lambda_1$ and $\Lambda_2$ are disjoint Borel
subsets of $\RR^d$, we obtain $a_{2,\infty}(m)=0$, and so (v) is
satisfied. Third, $\Sigma_n$ reduces to
\[
\Sigma_n=\int_{W_n}f_\theta^{(n)}(u)f_\theta^{(n)}(u)^\top
\rho(u)\dd u.
\]
\subsubsection{Main result}\label{secmainresult}

We now state our main 
result concerning the asymptotics for
the variational estimator based on $\bX_{W_n}$, that is, the estimator
%
\begin{equation}
\label{eqdefEst} \ParE= - A_{n}(\bX)^{-1} b_{n}(
\bX)
\end{equation}
defined when $A_{n}(\bX)=\widehat{S}_n$ given by
\[
\widehat{S}_n= \sum_{u \in\bX_{W_n}}
h^{(n)}(u) \div z(u)^\top
\]
is invertible, and where
\[
b_{n}(\bX)=\sum_{u \in\bX_{W_n}} \div
h^{(n)}(u).
\]
Denote $\xrightarrow{d}$ convergence in distribution as
$n\rightarrow\infty$. 

\begin{theorem}\label{thmconv} For $d\ge2$ and under the conditions \textup{(i)--(vi)},
the variational estimator $\ParE$ defined by (\ref{eqdefEst})
satisfies the following properties.

\textup{(a)} With probability one, when $n$ is sufficiently large,
$\widehat{S}_n$ is invertible
(and hence
$\ParE$ exists).

\textup{(b)}
$\ParE$ is a strongly consistent estimator of $\Par$.

\textup{(c)} We have 
%
\begin{equation}
\label{eqconvDist} {\Sigma}_n^{-1/2} {S}_n (
\ParE-\Par) \xrightarrow{d} \mathcal{N}(0,I_p),
\end{equation}
%
where ${\Sigma}_n^{-1/2}$ is the inverse of $\Sigma_n^{1/2}$, where
${\Sigma}_n^{1/2}$ is any square
matrix with
\mbox{${\Sigma}_n^{1/2}({\Sigma}_n^{1/2})^\top={\Sigma}_n$}.
\end{theorem}


Theorem \ref{thmconv} is verified in Appendix \ref{appA}, where,
for example, in the proof of Lemma \ref{lemZn} it becomes convenient that
$d\ge2$. We claim that the results of Theorem \ref{thmconv} remain
valid when $d=1$, but other conditions and another proof
are then needed, and we omit
these technical details.

\section{Simulation study}\label{secsim}

\subsection{Planar results with a modest number of points} \label{secsimd2}

In this section, we investigate the finite-sample properties of the
variational estimator (\textsc{vare}) for the planar case $d=2$ of an
inhomogeneous
Poisson point process, for an inhomogeneous log-Gaussian Cox process,
and for an
inhomogeneous Thomas process. We
compare \textsc{vare} with the maximum first-order composite likelihood estimator
(\textsc{mcle}) obtained by maximizing the composite log-likelihood
(discussed at the beginning of Section~\ref{secintro})
and which
is equivalent to the Poisson log-likelihood
%
\begin{equation}
\label{eqcle} \sum_{u \in \mathbf{X}_W} \log\rho(u) - \int
_{W} \rho(u) \dd u.
\end{equation}
In contrast to the variational approach, this provides not only an
estimator of $\Par$ but also of $\beta$.

It seems fair to compare the \textsc{vare} and the \textsc{mcle} since
both estimators are based only on the parametric
model for the log-linear intensity
function $\rho$. Guan and Shen \cite{guanshen10} and
Guan, Jalilian and Waagepetersen \cite{waagepetersenguanjalilian11}
show that the \textsc{mcle}
can be improved if a parametric model for the second order product
density $\rho^{(2)}$ is included when constructing a second-order
composite log-likelihood based
on both $\rho$ and $\rho^{(2)}$. We leave it as an open problem how to
improve our
variational approach by incorporating a parametric
model for $\rho^{(2)}$.

\begin{figure}

\includegraphics{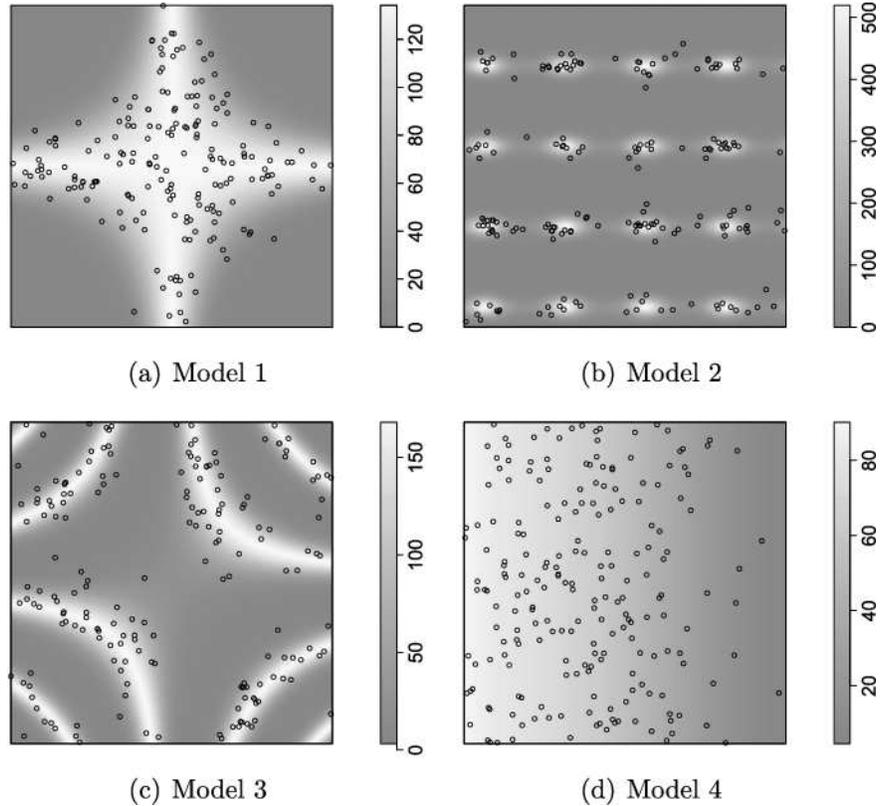}

\caption{Intensity functions and
examples of realizations of Poisson point
processes with intensity functions given by models 1--4 (defined in
Section \protect\ref{secsimd2}) and generated on the region $[-1,1]^2$.}
\label{figexModels}
\end{figure}

We consider four different models for the log-linear intensity
function given by (\ref{erho}), where $p=
1, 2, 1, 3$, respectively, and $u=(u_1,u_2)\in[-2,2]^2$:
\begin{itemize}
\item{Model 1}: $\Par=-2$, $z(u)=u_1^2u_2^2$.
\item{Model 2}: $\Par=(1,4)^\top$, $z(u)=( \sin(4\uppi  u_1) , \sin
(4\uppi
u_2) )^\top$.
\item{Model 3}: $\Par=2$, $z(u)= \sin(4\uppi  u_1u_2)$.
\item{Model 4}: $\Par=(-1,-1,-0.5)^\top$,
$z(u)=(u_1,u_1^2,u_1^3)^\top$.
\end{itemize}
%
We assume that the covariate function $z(u)$ is known to us for all
$u\in W$ so that
we can evaluate its first and second derivatives (Section \ref{secz-finite}
considers the case where $z$ is only known at a finite set of locations).
Figure \ref{figexModels} shows the intensity functions and
simulated point patterns under
models 1--4 for a Poisson point process within the region $W=[-1,1]^2$.
The figure illustrates the
different types of
inhomogeneity obtained by the different choices of $\rho$. 

In addition to the Poisson point process, referred to as \textsc{poisson}
in the
results to follow, two cases of
Cox process models
are considered, where we are using the terminology and notation introduced
in Section \ref{secdisccond}:
\begin{itemize}
\item An inhomogeneous log-Gaussian Cox process $\bX$ where the
underlying Gaussian process has an exponential covariance function
$c(u,v)=\sigma^2\exp(-\|u-v\|/\alpha)$. We
refer then to $\bX$
as \textsc{lgcp1}
when $\sigma^2=0.5$ and $\alpha=1/15$, and as \textsc{lgcp2} when
$\sigma^2=1.5$ and $\alpha=1/30$.
%
\item An inhomogeneous Thomas process $\bX$ where $\kappa$ is the
intensity of the underlying Poisson point process $\bC$ and
$\sigma$ is the standard deviation of the normal density
$f_\sigma$, see (\ref{einhomSNCP}). We refer then to $\bX$
as \textsc{thomas1} when $\kappa=100$ and
$\sigma=0.05$, and as \textsc{thomas2} when $\kappa=300$ and $\sigma=0.1$.
\end{itemize}
%

In addition two observation windows are considered: $W=W_1=[-1,1]^2$ and
$W=W_2=[-2,2]^2$. For each choice of model and observation window,
we adjusted the parameter $\beta$
such that the expected number of points, denoted by
$\mu^\star$, is 200 for the choice $W=W_1$ and 800 for the choice $W=W_2$
(reflecting the fact that $W_2$ is four times larger than $W_1$),
and then
1000 independent point patterns were
simulated using the \texttt{spatstat} package of \texttt{R}
Baddeley and Turner \cite{baddeleyturner05}.

For each of such 1000 replications, we computed
the \textsc{mcle}, using the \texttt{ppm()} function of
\texttt{spatstat} with a fixed deterministic grid of $80 \times80$
points to discretize the integral in (\ref{eqcle}).
We also
computed
the \textsc{vare} considering either
the test function $h(u)=\div z(u)$
or its modification $h(u)=\div z(u) \eta_W(u)$ for various values of
$\varepsilon>0$, where the
former case can be viewed as a limiting case of the latter one
with $\varepsilon=0$. For the other choices of test functions
discussed in
Section \ref{secnewsec}
some preliminary experiments showed that the present choice
of test functions led to estimators with the smallest variances.

Among the different models for the intensity function, models 2 and 4
are indeed correctly defined on $\RR^d$
in the sense that they satisfy at least our condition (ii). To
illustrate the simplicity of
the \textsc{vare} and the flexibility of conditions (i)--(vi), we focus on
model 2 in Appendix \ref{appmodel2}, detail
the form of the \textsc{vare}, and show that our asymptotic results are valid.

\begin{figure}

\includegraphics{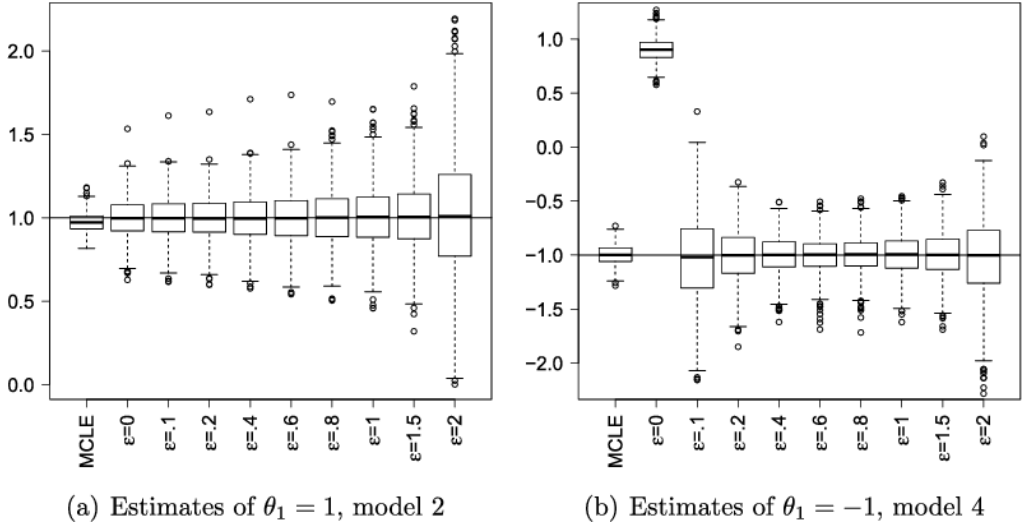}

\caption{Box plots of \textsc{mcle} and \textsc{vare} for
$\theta_1$ (the first coordinate of $\theta$) under models 2 and 4, when
using the test function $h(u)=\div z(u) \eta_W(u)$ for
different values of $\varepsilon$, with $\varepsilon=0$ corresponding
to $h(u)=\div z(u)$. The plots are based on simulations from
Poisson point processes on the observation window $[-2,2]^2$, when the expected
number of points is 800. Similar results are obtained for the other
cases of point process models and choice of observation window.}
\label{figbpl}
\end{figure}

Figure \ref{figbpl} illustrates some general findings for any choice
of point process model and observation window:
When the smoothing parameter $\varepsilon$ is at least $5\%$ larger
than the side-length of
the observation window, the \textsc{vare} is effectively
unbiased, and its variance increases as $\varepsilon$ increases.
However, when the point process is too much aggregated on the boundary
of the
observation window (as, e.g., in the case of
(b) in Figure \ref{figexModels}), a too
small value of $\varepsilon$ leads to biased estimates. At the
opposite, when the point process is not too much aggregated on the
boundary of the observation window (see, e.g., in the case of (a) in
Figure \ref{figexModels}), the choice $\varepsilon=0$ leads to
the smallest variance.

\begin{table}
\caption{Average of the
$p$ empirical mean squared
errors (\textsc{amse})
of the estimates for the coordinates in
$\Par=(\theta_1,\ldots,\theta_p)^\top$ and
based on independent realizations of Poisson, inhomogeneous
log-Gaussian Cox processes, and
inhomogeneous Thomas point processes with different parameters,
intensity functions, and observation windows as described in
Section \protect\ref{secsimd2}} 
\label{tabresd2}
\begin{tabular*}{\tablewidth}{@{\extracolsep{4in minus 4in}}lllllll@{}}
\hline
&\multicolumn{3}{l}{$W_1=[-1,1]^2$ $(\mu^\star=200)$}
&\multicolumn{3}{l@{}}{$W_2 = [-2,2]^2$ $(\mu^\star=800)$} \\[-4pt]
&\multicolumn{3}{c}{\hrulefill}
&\multicolumn{3}{c@{}}{\hrulefill} \\
%
%
& \multicolumn{2}{l}{\textsc{vare}} &  &
\multicolumn{2}{l}{\textsc{vare}} &  \\[-4pt]
& \multicolumn{2}{c}{\hrulefill} &  &
\multicolumn{2}{c}{\hrulefill} &
\\
& $\varepsilon=0$ & $\varepsilon=0.1$ & \textsc{mcle} &
$\varepsilon=0$ & $\varepsilon=0.2$ & \textsc{mcle}\\
\hline
\multicolumn{7}{@{}c@{}}{\textit{Model} 1: $\Par=-2$,
$z(u)=u_1^2u_2^2$}\\[4pt]
\textsc{poisson} & 0.109 & 0.124 & 0.085 & 0.027 &
0.030 & 0.022 \\
\textsc{lgcp1} & 0.152 & 0.181 & 0.143 & 0.035 & 0.040 & 0.032 \\
\textsc{lgcp2} & 0.170 & 0.203 & 0.143 & 0.035 & 0.041 & 0.033 \\
\textsc{thomas1} & 0.141 & 0.163 & 0.118 & 0.033 & 0.037 & 0.030 \\
\textsc{thomas2} & 0.118 & 0.147 & 0.095 & 0.026 & 0.027 & 0.025 \\
[6pt]
\multicolumn{7}{@{}c@{}}{\textit{Model} 2: $\Par=(1,4)^\top$, $z(u)=( \sin
(4\uppi u_1) , \sin(4\uppi  u_2) )^\top$} \\[4pt]
\textsc{poisson} & 0.104 & 0.126 & 0.089 & 0.028 & 0.033 & 0.033 \\
\textsc{lgcp1}& 0.131 & 0.159 & 0.117 & 0.041 & 0.047 & 0.066 \\
\textsc{lgcp2} & 0.180 & 0.213 & 0.144 & 0.055 & 0.062 & 0.067 \\
\textsc{thomas1} & 0.132 & 0.158 & 0.106 & 0.039 & 0.046 & 0.062 \\
\textsc{thomas2} & 0.106 & 0.130 & 0.098 & 0.035 & 0.039 & 0.061 \\
[6pt]
\multicolumn{7}{@{}c@{}}{\textit{Model} 3: $\Par=2$, $z(u)= \sin(4\uppi  u_1u_2)$}\\
[4pt]
\textsc{poisson} & 0.087 & 0.105 & 0.037 & 0.023 & 0.026 & 0.010 \\
\textsc{lgcp1} & 0.122 & 0.137 & 0.052 & 0.038 & 0.036 & 0.023 \\
\textsc{lgcp2} & 0.149 & 0.174 & 0.057 & 0.038 & 0.038 & 0.023 \\
\textsc{thomas1} & 0.103 & 0.119 & 0.048 & 0.033 & 0.032 & 0.021 \\
\textsc{thomas2} & 0.096 & 0.109 & 0.042 & 0.034 & 0.031 & 0.021 \\
[6pt]
\multicolumn{7}{@{}c@{}}{\textit{Model} 4: $\Par=(-1,-1,-0.5)^\top$,
$z(u)=(u_1,u_1^2,u_1^3)^\top$}\\[4pt]
\textsc{poisson} & 0.420 & 0.410 & 0.216 & 1.819 & 0.027 & 0.010 \\
\textsc{lgcp1} & 0.463 & 0.556 & 0.332 & 1.835 & 0.035 & 0.015 \\
\textsc{lgcp2} & 0.471 & 0.588 & 0.327 & 1.841 & 0.035 & 0.016 \\
\textsc{thomas1} & 0.456 & 0.545 & 0.277 & 1.836 & 0.030 & 0.012 \\
\textsc{thomas2} & 0.427 & 0.445 & 0.246 & 1.805 & 0.026 & 0.010 \\
\hline
\end{tabular*}
\end{table}

Table \ref{tabresd2} concerns the situations with
$\varepsilon=0$, $\varepsilon=0.1$ when $W=W_1=[-1,1]^2$, and
$\varepsilon=0.2$ when $W=W_2=[-2,2]^2$ (in the latter two cases, the
choice of
$\varepsilon>0$ corresponds to $5\%$ of the side-length of $W$). The table
shows the average of the $p$ empirical mean\vadjust{\goodbreak} squared
errors (abbreviated as \textsc{amse})
of the estimates for the coordinates in $\Par=(\theta_1,\ldots,\theta
_p)^\top$
and based on the 1000 replications. In all except a few cases, the
\textsc{amse} is smallest for the \textsc{mcle}, the exception being model 2
when $W=W_2$. In most cases, the
\textsc{amse} is smaller when $\varepsilon=0$ than if $\varepsilon>0$,
the exception being some cases of model 3 when $W=W_2$ and all cases
of model 4 when $W=W_2$. For models 1--2, the \textsc{amse} for the \textsc{vare} with
$\varepsilon=0$ is rather close to the \textsc{amse} for the \textsc{mcle}. For
models 3--4, and in particular model 4 with $W=W_2$, the difference is
more pronounced, and the \textsc{amse} for the \textsc{mcle} is the smallest.

\subsection{Results with a high number of points and varying
dimension of space}\label{secres2}

In this section, we investigate the \textsc{vare} and the \textsc{mcle} when
the observed number of points is expected to be very high, when
the dimension $d$ varies from 2 to 6, and when
the dimension $p$ of $\theta$ scales with $d$.
Specifically, we let $p=d$ and consider a Poisson point process with
\[
\log\rho(u) =\beta+ \sum_{i=1}^d
\theta_i\sin(4\uppi  u_i)/d,\qquad u=(u_1,\ldots,u_d)^\top\in\RR^d,
\]
where $\theta_1=\cdots=\theta_d=1$, $d=2,\ldots,6$,
and
$\beta$ is chosen such that
the expected number of points in $W=[-1,1]^d$ is
$\mu^\star=10\,000$.

For $d=2,\ldots,6$,
we simulated 1000 independent realizations of such a Poisson point
process within
$W=[-1,1]^d$. For each realization,
when
calculating the \textsc{mcle}
we used a systematic
grid (i.e., a square, cubic$,\ldots$ grid when $d=2,3,\ldots$)
for the discretization of the integral in~(\ref{eqcle}),
where the number of dummy points $n_D$ is equal to $\tau\mu^\star$ with
$\tau=0.1,0.5,1,2,4,10$.

Similar to Table \ref{tabresd2},
Table \ref{tabmsedimd} shows ratios of \textsc{amse}'s for the two
types of estimators, \textsc{vare} and \textsc{mcle}, as the dimension $d$
(and number of parameters) varies and as
the number of dummy points $n_D$ varies from 1000 to
100\,000. In terms of the
\textsc{amse}, the \textsc{vare} outperforms the \textsc{mcle} for the smaller
values of $n_D$, and the two estimators are only equally good at the
largest value of $n_D$ in Table \ref{tabmsedimd}.


\begin{table}[b]
\caption{Ratio of the \textsc{amse} of the \textsc{mcle} over the \textsc{amse}
of the \textsc{vare} for 
$\theta=(\theta_1,\ldots,\theta_d)\in\RR^d$ and based on
simulations from
Poisson point processes as described in Section \protect\ref{secres2}.
The rows corresponds to the dimension (and number of
parameters) $d$, and the columns to
the number of dummy points $n_D=10\,000\tau$ used to
discretize the integral of
(\protect\ref{eqcle}) when calculating the \textsc{mcle}}
\label{tabmsedimd}
\begin{tabular*}{\tablewidth}{@{\extracolsep{\fill}}lllllll@{}}
\hline
& \multicolumn{6}{l@{}}{$\mbox{\textsc{amse}}_{\mathrm{MCLE}}/
\mbox{\textsc{amse}}_{\mathrm{VARE}}$} 
\\[-4pt]
& \multicolumn{6}{l@{}}{\hrulefill}\\
& $\tau=0.1$ & $\tau=0.5$ & $\tau=1$ & $\tau=2$ & $\tau=4$ &$\tau=10$
\\
\hline
$d=2$ & 11.00 & 2.71 & 1.83 & 1.32 & 1.08 & 0.95 \\
$d=3$ & 11.20 & 2.77 & 1.88 & 1.36 & 1.15 & 0.99 \\
$d=4$ & 11.35 & 2.92 & 1.97 & 1.41 & 1.16 & 0.99 \\
$d=5$ & 11.67 & 3.00 & 2.00 & 1.43 & 1.21 & 1.03 \\
$d=6$ & 10.59 & 2.92 & 1.92 & 1.40 & 1.17 & 1.02 \\
\hline
\end{tabular*}
\end{table}

\begin{table}
\caption{Average time (in seconds) for the
computation of the \textsc{vare} and of the \textsc{mcle} as considered in
Table \protect\ref{tabmsedimd}} 
\label{tabtime}
\begin{tabular*}{\tablewidth}{@{\extracolsep{\fill}}llllllll@{}}
\hline
& & \multicolumn{6}{l@{}}{\textsc{mcle}} \\[-4pt]
& & \multicolumn{6}{l@{}}{\hrulefill}
\\
&\textsc{vare} & $\tau=0.1$ & $\tau=0.5$ & $\tau=1$ & $\tau=2$ & $\tau=4$ &$\tau=10$
\\
\hline
$d=2$ & 0.004 & 0.200 & 0.347 & 0.546 & 0.984 & 1.929 & 5.744 \\
$d=3$ & 0.005 & 0.178 & 0.298 & 0.450 & 0.779 & 1.483 & 4.087 \\
$d=4$ & 0.007 & 0.231 & 0.374 & 0.562 & 0.941 & 1.740 & 4.805 \\
$d=5$ & 0.009 & 0.272 & 0.432 & 0.650 & 1.082 & 1.994 & 5.493 \\
$d=6$ & 0.011 & 0.312 & 0.494 & 0.739 & 1.242 & 2.367 & 6.203 \\
\hline
\end{tabular*}
\end{table}

Table \ref{tabtime} presents the average time in seconds to get one
estimate based on the \textsc{vare} and as a function of $d$, and also
the average time in seconds to get one
estimate based on the \textsc{mcle} and
as a function of both $d$ and $\tau$. The table clearly shows how much
faster the calculation of the \textsc{vare} than the \textsc{mcle} is.
In particular, when
$n_D=100\,000$, the average computation time of the \textsc{mcle} is
around 1400 ($d=2$) to 560 ($d=6$) times slower than that of the \textsc{vare}.

\subsection{Results when $z$ is known only on a finite set of
locations}\label{secz-finite}

The calculation of the \textsc{vare} based on a realization $\bX
_W=\mathbf{x}$
requires the knowledge of $\div z(u)$ (and possibly also $\div\div z(u)$)
for $u\in\mathbf{x}$. In practice, $z$ is often
only known for a finite set of
points in $W$, which is usually given by a systematic grid
imposed on $W$, and we propose then to
approximate $\div z$ and
$\div\div z$ using the finite-difference method. We discuss below
some interesting findings when such an approximation is used.

\begin{figure}

\includegraphics{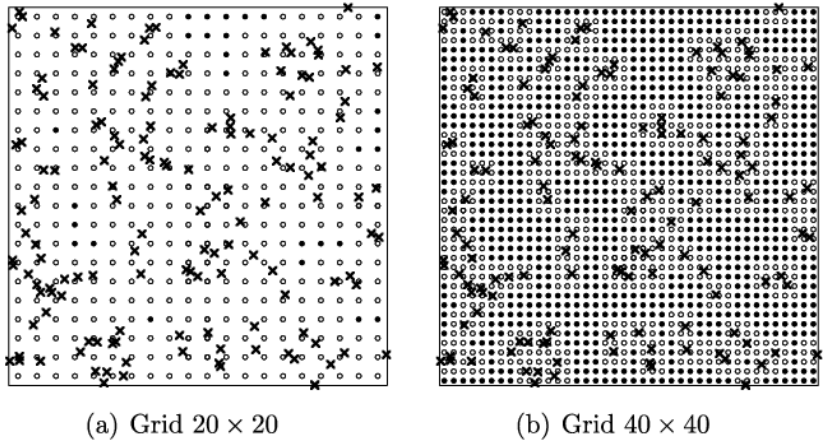}

\caption{The crosses represent a realization of the Poisson point
process under the model 2 and within the observation window $[-1,1]^2$.
The empty circles represent the grid points where the spatial function
$z$ is sampled and used to compute \textsc{vare}(loc) and \textsc{mcle}(loc).
The grid points used to compute the \textsc{mcle} correspond to the empty
and filled circles.} \label{figexCov}
\end{figure}

We focus on the planar case $d=2$, and let $h(u)=\div z(u)$ for the
\textsc{vare}. For the
two choices of observation windows, $W=W_1=[-1,1]^2$ or $W=W_2=[-2,2]^2$,
we simulated 1000 realizations of a
Poisson point process with $\log\rho(u)=\beta+\sin(4\uppi u_1)+\sin
(4\uppi
u_2)$ for
$u=(u_1,u_2)\in\RR^2$
(i.e., model 2 in Section \ref{secsimd2} with $\theta_1=\theta_2=1$),
where $\beta$ is
chosen such that the expected number of points is $\mu^\star=200$ if
$W=W_1$ and $\mu^\star=800$ if
$W=W_2$. For each replication, we calculated four types
of
estimators, namely \textsc{vare} and \textsc{mcle}
which correspond to the situation in
Table \ref{tabresd2} where $z$ is assumed to be known on $W$,
and two ``local'' versions
\textsc{vare}(loc) and \textsc{mcle}(loc) where only
knowledge about $z$ on a grid is used. In detail:
\begin{itemize}
\item Assuming the full information about $z$ on $W$, \textsc{vare} and
\textsc{mcle} were calculated, where for the \textsc{mcle} the integral
in (\ref{eqcle}) is discretized over a quadratic grid $G$ of $n_D^2$
points in $W$, with $n_D=20,40,80$ if $W=W_1$, and $n_D=40,80,160$ if
$W=W_2$.
\item
For
each simulated point $u$ of a replication, the $3\times3$ subgrid whose
midpoint is closest to $u$ was used for
approximating $\div z(u)$ and
$\div\div z(u)$ by the finite-difference method. Thereby, a subgrid
$G_0\subseteq
G$ was obtained as illustrated in
Figure \ref{figexCov}. Using only the knowledge
about $z$ on $G_0$, \textsc{vare}(loc)
as an approximation of \textsc{vare} was obtained. Furthermore,
\textsc{mcle}(loc) was calculated by discretizing
the integral in (\ref{eqcle}) over the grid points in $G_0$.
\end{itemize}


\begin{table}[b]
\caption{\textsc{amse} for the four types of estimators
\textsc{vare}, \textsc{vare}(loc), \textsc{mcle}, and \textsc{mcle}(loc) obtained
using different grids
as described in Section \protect\ref{secz-finite}.
The \textsc{vare} is assuming that the spatial function $z$ is known and
is used here as a reference; it does not depend on the refinement of
the grid.
The results are based on 1000 independent\vspace*{1pt} realizations of a
planar Poisson point process simulated on the observation
window $W=[-1,1]^2$ or
$W=[-2,2]^2$}
\label{tabgrid}
\begin{tabular*}{\tablewidth}{@{\extracolsep{\fill}}lllllll@{}}
\hline
& \multicolumn{3}{l}{$W=[-1,1]^2$ $(\mu^\star=200)$}
&\multicolumn{3}{l@{}}{$W=[-2,2]^2$ $(\mu^\star=800)$} \\[-4pt]
& \multicolumn{3}{l}{\hrulefill} & \multicolumn{3}{l@{}}{\hrulefill}
\\
& $20 \times20$ &$40\times40$&$80\times80$ & $40\times40$
&$80\times80$&$160\times160$\\
\hline
\textsc{vare} & $-$ &{0.023} & $-$ & $-$&{0.006} &$-$\\ \textsc{vare}(loc) &
0.072 & 0.029 & 0.025 & 0.035 & 0.008 & 0.006 \\ \textsc{mcle} & 0.014 &
0.014 & 0.013 & 0.004 & 0.004 & 0.003 \\ \textsc{mcle}(loc) & 0.014 &
0.166 & 0.628 & 0.004 & 0.164 & 0.623 \\
\hline
\end{tabular*}
\end{table}

Table \ref{tabgrid} shows that in terms of the \textsc{amse}, the \textsc{vare}(loc)
is effectively as good as the \textsc{vare} if the grid is sufficiently
fine, cf. the results in the case of the $80\times80$ grid for $W_1$
and the $160\times160$ grid for $W_2$. As expected the \textsc{mcle}
performs better than the other estimators, in particular as the grid
becomes finer, except for the coarsest grids\vadjust{\goodbreak} (the~$20\times20$ grid
for $W_1$
and the $40\times40$ grid for~$W_2$) where the \textsc{amse} is equal for the
\textsc{mcle} and the \textsc{mcle}(loc). As the grid gets finer,
the \textsc{amse} for the \textsc{mcle}(loc) increases and becomes much larger
than for any of the other estimators -- only for the coarsest grids,
the \textsc{mcle}(loc) and the \textsc{mcle} perform equally good. Thus
if the covariates are observed only in a small neighborhood of
the location points, it becomes advantageous to use
the \textsc{vare} as compared to the \textsc{mcle}.
This feature could be of relevance in practice if the covariates are
only determined at locations close to the points of $\bX_W$.

\begin{appendix}\label{app}
\section{Proofs}\label{appA}

This Appendix verifies Theorem \ref{thmconv} and some accompanying
lemmas assuming that $d\ge2$ and
conditions (i)--(vi) in Section \ref{secconditions}
are satisfied.

To simplify the notation, when
considering a mean value which possibly
depends on $(\beta,\Par)$, we suppress this and simply write
$\E[\cdots]$.

We start by showing that we can replace:
\begin{enumerate}
\item the domain $W_n$ by a more convenient domain $W_n^\star$
satisfying $|W_n|\sim|W_n^\star|$ as $n\to\infty$ (meaning that
$|W_n|/|W_n^\star|\rightarrow1$ as $n\to\infty$);
\item the function $h^{(n)}$ by a function $h_\varepsilon^{(n)}$ with
compact support on $W_n^\star$, where $\varepsilon=\varepsilon_n$
depends on $n$ and should be distinct from the $\varepsilon$ used in
(\ref{eeta}).
\end{enumerate}
This will later allow us to apply Corollary \ref{corest}.

Let $C_i=i+(-1/2,1/2]^d$ be the unit box centered at
$i\in\ZZ^d$. Define $\mathcal I_n=\{i\in\ZZ^d\dvt
C_i\subset W_n\}$, and let $\partial\mathcal I_n =\{
i\in\ZZ^d\setminus\mathcal I_n\dvt  C_i\cap W_n\neq\varnothing\}$ be the
nearest neighbourhood of $\mathcal I_n$ on
the integer lattice $\ZZ^d$. Set $W_{\partial\mathcal
I_n}=\bigcup_{i\in\partial\mathcal I_n}C_i$ and $W_n^\star=
\bigcup_{i\in
\mathcal I_n} C_i$.

\begin{lemma} \label{lemdomain} For any $n=1,2,\ldots\,$, we have
$W_n^\star\subseteq W_n
\subseteq W_n^\star\cup W_{\partial\mathcal I_n}$.
As $n\to\infty$, then $|W_n|= |A| n^d\sim|W_n^\star|$ and
$|W_n\setminus W_n^\star| =\mathcal O(n^{d-1})$. Moreover, $\sum_{n\geq
1}|\mathcal I_n|^{-1}<\infty$.
\end{lemma}

\begin{pf} The first statement is clearly true. Thus,\vspace*{1pt}
$|W_n^\star|\leq|W_n|\leq|W_n^\star|+|\partial\mathcal I_n|$.

By (i), $W_n=nA$ is convex, so
$|\partial\mathcal I_n| \leq K_d\delta(A)^{d-1} n^{d-1}$,
where $\delta(A)$ denotes the diameter of $A$ and $K_d>0$ is a
constant. Consequently,
\[
1 \geq\frac{|W_n^\star|}{|W_n|} \geq1 - \frac{|\partial\mathcal
I_n|}{|W_n|} \geq1-
\frac{\delta(A)^{d-1}}{n}
\]
leading to\vspace*{1pt} $|W_n|\sim|W_n^\star|$ as $n\to\infty$.
Since $|W_n\setminus W_n^\star|/|W_n| \leq K_d\delta(A)^{d-1}/n =
\mathcal{O}(1/n)$, we obtain $|W_n\setminus W_n^\star|=\mathcal
O(n^{d-1})$, whereby the second statement is verified.

The last statement follows from that $|\mathcal
I_n|=|W_n^\star|\sim|A|n^d$ and $d\ge2$.
\end{pf}

Now, let $\varepsilon=\varepsilon_n=n^\alpha$ for some given $\alpha
\in[0,1)$. Define $h_\varepsilon^{(n)}$ as the regularized function
of $h^{(n)}$ as described in Section \ref{secnew-est} and given by
%
\begin{equation}
\label{eqhepsilon} h_\varepsilon^{(n)} (u) = h^{(n)}(u)
\eta_{W_n^\star}(u),
\end{equation}
where $\eta_{W_n^\star}$ is defined by (\ref{eeta}) (when $W$ is
replaced by $W_n^\star$ and the $\varepsilon$ in (\ref{eeta}) is replaced
by the present $\varepsilon=\varepsilon_n$).
By Lemma \ref{lemhepsilon} and (i)--(ii), we have that
$h_\varepsilon^{(n)}$ respective $\div h_\varepsilon^{(n)}$ agrees
with $h^{(n)}$ respective $\div h^{(n)}$ on ${W_n^\star}_{\ominus
2\varepsilon}$, the support of $h_\varepsilon^{(n)}$ is included in the
bounded set $W_n^\star$, and there exists $K<\infty$ such that
%
\begin{equation}
\label{eqhen} \sup_{n\geq1} \bigl\|h_\varepsilon^{(n)}
\bigr\|_\infty\leq K \quad\mbox{and}\quad \sup_{n\geq1} \bigl\|\div
h_\varepsilon^{(n)} - \div h^{(n)} \bigr\|_\infty\leq
K.
\end{equation}

The following lemma concerns the behavior of variance functionals
computed on $W_n$ or $W_n^\star$.

\begin{lemma}\label{lemvar}
Let $(\psi^{(n)})_{n\geq1}$ be a sequence of functions in $\mathcal
C^0_{d,1}$ such that
%
\begin{equation}
\label{eqpsin} \sup_{n\geq1} \bigl\| \psi^{(n)}
\bigr\|_\infty\leq C
\end{equation}
for some constant $C<\infty$, then for $\widetilde W_n=W_n,W_n^\star$,
the variance
\[
V_{\widetilde W_n} = \Var\biggl( \sum_{u \in\bX_{\widetilde W_n}}
\psi^{(n)}(u) \biggr)
\]
is finite and is given by
%
\begin{equation}
\label{eqvar} V_{\widetilde W_n} 
= \int_{\widetilde W_n}
\psi^{(n)}(u)^2 \rho(u) \dd u + \int_{\widetilde W_n}
\int_{\widetilde W_n} \psi^{(n)}(u) \psi^{(n)}(v)
Q_2(u,v)\dd u \dd v = \mathcal O\bigl(n^d\bigr).
\end{equation}
\end{lemma}

\begin{pf} The finiteness of the variance follows from (iv), and
the first identity in (\ref{eqvar}) is immediately derived from
(\ref{ecampbell})--(\ref{ecampbell2}).

For the second identity,
we consider first $\widetilde W_n=W_n^\star$. Define
$Y_i^{(n)} = \sum_{u \in C_i} \psi^{(n)}(u)$ for $i \in\mathcal
I_n$. For $\delta\geq1$ given in (iv), it is clear that $\E
(|Y_i^{(n)}|^{2+\delta})$ is bounded by a linear combination of
\[
s_k^{(n)}= \int_{C_i}\cdots\int _{C_i} \bigl|\psi^{(n)}(u_1)\cdots
\psi^{(n)}(u_{k})\bigr|\rho^{(k)}(u_1,\ldots,u_k)\dd u_1
\cdots\mathrm{d} u_k,\qquad k=1,\ldots,2+\delta.
\]
Using (\ref{eqpsin}) and (iv), we obtain
\[
\sup_{n\geq1} s_k^{(n)} \leq C^k
\sup_{i \in\ZZ^d} \int_{C_i}\cdots\int
_{C_i} \rho^{(k)}(u_1,\ldots,u_k)\dd u_1 \cdots\mathrm{d} u_k \le
C^kK'<\infty.
\]
Therefore,
\[
M_Y:=\sup_{n\geq1}\sup_{i \in\mathcal I_n} \E
\bigl(\bigl|Y_i^{(n)}\bigr|^{2+\delta} \bigr) <\infty.
\]
Further, we have the following bound for the covariance in terms of the
mixing coefficients of $\bX$ (see Doukhan \cite{doukhan94} or Guyon
\cite{guyon91}, remark, page 110),
\[
\bigl|\Cov\bigl( Y_i^{(n)}, Y_j^{(n)}
\bigr) \bigr| \leq8M_Y^2 \alpha_{1,1}\bigl(|j-i|\bigr)^{{\delta}/({2+\delta})}.
\]
Furthermore, since for any $m\geq1$, $\alpha_{1,1}(m)\leq\alpha
_{2,\infty}(m)$, and since $|W_n^\star|=|\mathcal I_n|$, we obtain
\begin{eqnarray*}
\bigl|W_n^\star\bigr|^{-1} V_{W_n^\star} &=& |\mathcal
I_n|^{-1} \sum_{i,j\in
\mathcal I_n} \Cov
\bigl(Y_i^{(n)},Y_j^{(n)}\bigr)
\\
&\leq&8M_Y^2 |\mathcal I_n|^{-1}
\sum_{i,j \in\mathcal I_n} \alpha_{2,\infty}\bigl(|j-i|\bigr)^{
\delta/(2+\delta)}
\\
&\leq&8M_Y^2 \sum_{m\geq0} \bigl|
\bigl\{j \in\ZZ^d\dvt  |j|=m \bigr\} \bigr| \alpha_{2,\infty}
(m)^{\delta/(2+\delta)}
\\
&\leq& c_d \sum_{m\geq1} m^{d-1}
\alpha_{2,\infty} (m)^{\delta/(2+\delta)},
\end{eqnarray*}
where $c_d>0$ is a constant depending only on $d$. Combining this with
(v) leads to $|W_n^\star|^{-1} V_{W_n^\star}=\mathcal O(1)$.

Second, let $\mathcal J_n= \mathcal I_n \cup\partial\mathcal I_n$. Then
\[
V_{W_n} = \sum_{i,j\in\mathcal J_n} \Cov\bigl(
Z_i^{(n)}, Z_j^{(n)} \bigr) \qquad\mbox{where for } i\in\mathcal J_n,\qquad
Z_i^{(n)} = \sum
_{u\in X_{C_i\cap W_n} } \psi^{(n)}(u).
\]
Using (\ref{eqpsin}), (iv), and similar arguments as above for
the case $\widetilde W_n=W_n$, it is clear that
\[
M_Z:=\sup_{n\geq1}\sup_{i \in\mathcal J_n} \E
\bigl(\bigl|Z_i^{(n)}\bigr|^{2+\delta} \bigr) <\infty.
\]
Finally, using (v) and similar arguments as above, we obtain that
$|\mathcal
J_n|^{-1}V_{W_n}=\mathcal O(1)$. This completes the proof, since
$|\mathcal J_n|\sim|\mathcal I_n|=\mathcal O( n^d)$.
\end{pf}

Similar to the definitions of $A_{n}(\bX)$ and $b_{n}(\bX)$ in
Section \ref{secnew-est}, we define
\[
A_{n}^\star(\bX) = \sum_{u \in\bX_{W_n^\star}}
h_\varepsilon^{(n)}(u)\div z(u)^\top\quad\mbox{and}\quad
b_{n}^\star(\bX) = \sum_{u \in\bX_{W_n^\star}}
\div h_\varepsilon^{(n)}(u).
\]
%
We simplify the notation by suppressing the dependence on
$\bX$ for the random matrices $A_{n}=A_{n}(\bX)$ and
$A_{n}^\star=A_{n}^\star(\bX)$,
and for the random vectors $b_{n}=b_{n}(\bX)$ and
$b_{n}^\star=b_{n}^\star(\bX)$. 

\begin{lemma}\label{lemZn}
\textup{(I)} For $Z_n=A_{n},A_{n}^\star,b_{n},b_{n}^\star$, we have
$|W_n|^{-1}(Z_n-\E Z_n) \xrightarrow{\mathit{a.s.}} 0$ as $n \to\infty$.

\textup{\hphantom{I}(II)} $|W_n|^{-1}\E( A_{n} \Par+ b_{n} ) = \mathcal O
({n^{\alpha-1}})$.

\textup{(III)} $ (A_{n}-A_{n}^\star)\Par+b_{n}-b_{n}^\star=
\RMo_{P}(|W_n|^{1/2})= \RMo_P(n^{d/2})$.
\end{lemma}

\begin{pf}
(I): We have
\begin{eqnarray*}
A_{n}-\E A_{n} &=& \biggl(\sum
_{u \in X_{W_n}} h^{(n)}(u) \div z(u)^\top\biggr)-
\int_{W_n}h^{(n)}(u)\div z(u)^\top\rho(u)
\dd u,
\\
A_{n}^\star-\E A_{n}^\star &=& \biggl(\sum
_{u \in X_{W_n^\star}} h_\varepsilon^{(n)}(u) \div
z(u)^\top\biggr) -\int_{W_n^\star
}h_\varepsilon^{(n)}(u)
\div z(u)^\top\rho(u) \dd u,
\\
b_{n}-\E b_{n} &=& \biggl(\sum
_{u \in X_{W_n}} \div h^{(n)} (u) \biggr) -\int
_{W_n}\div h^{(n)} (u) \rho(u) \dd u,
\\
b_{n}^\star-\E b_{n}^\star &=& \biggl(
\sum_{u \in X_{W_n^\star}}\div h_\varepsilon^{(n)}(u)
\biggr) -\int_{W_n^\star}\div h_\varepsilon^{(n)}(u)
\rho(u) \dd u .
\end{eqnarray*}
Let $j,k\in\{1,\ldots,p\}$.
From (ii) and (\ref{eqvar}), we obtain
\[
\E\bigl( (A_{n}-\E A_{n} )_{jk}^2
\bigr) = \mathcal O \bigl(n^d\bigr),\qquad \E\bigl( (b_{n}-\E
b_{n} )_{j}^2 \bigr)= \mathcal O
\bigl(n^d\bigr),
\]
%
%
\[
\E\bigl( \bigl(A_{n}^\star-\E
A_n^\star\bigr)_{jk}^2 \bigr) =
\mathcal O\bigl(n^d\bigr),\qquad 
\E\bigl( \bigl(b_n^\star-
\E b_n^\star\bigr)_{j}^2 \bigr)=
\mathcal O \bigl(n^d\bigr).
\]
Hence, for $Z_n=A_{n},A_n^\star,b_{n},b_n^\star$, we have (setting
$k=1$ for $Z_n=b_{n},b_n^\star$)
\[
\Var\bigl(|W_n|^{-1} ({Z_n})_{jk}
\bigr) = \mathcal O \bigl(n^{-d}\bigr),
\]
which together with the Borel--Cantelli lemma and the fact that $d\geq2$
imply the result of (I).

\mbox{}\hphantom{I}(II): By Lemma \ref{lemhepsilon} and
(\ref{eqhepsilon})--(\ref{eqhen}), we have 
%
\begin{equation}
\label{eqdiffA} A_{n}-A_n^\star= \sum
_{u \in\bX_{W_n^\star\setminus{W_n^\star
}_{\ominus2\varepsilon}}} \bigl(h^{(n)} (u)-h_\varepsilon^{(n)}(u)
\bigr) \div z(u)^\top+\sum_{u \in\bX_{W_n\setminus{W_n^\star}}}
h^{(n)}(u) \div z(u)^\top
\end{equation}
and 
%
\begin{equation}
\label{eqdiffb} b_{n}-b_{n}^\star= \sum
_{u \in
W_n^\star\setminus{W_n^\star}_{\ominus2\varepsilon}} \bigl(\div h^{(n)}
(u)-\div
h_\varepsilon^{(n)}(u)\bigr) +\sum_{u \in W_n\setminus{W_n^\star}}
\div h^{(n)} (u).
\end{equation}
We denote by $T_1$ and $T_2$ the two sums of the right-hand side
of (\ref{eqdiffA}) and by $T^\prime_1$ and $T^\prime_2$ the two sums
of the right-hand side of (\ref{eqdiffb}).
Using (ii), (\ref{ecampbell}), and (\ref{eqhen}), we obtain $\E T_1
= \mathcal O ( |W_n^\star\setminus{W_n^\star}_{\ominus
2\varepsilon}|)$, $\E T_2 = \mathcal O ( |W_n\setminus
{W_n^\star}|)$, $\E T_1^\prime= \mathcal O ( |W_n^\star\setminus
{W_n^\star}_{\ominus2\varepsilon}|)$, and $\E T_2^\prime= \mathcal
O ( |W_n\setminus{W_n^\star}|)$. By Lem\-ma~\ref{lemdomain},
$|W_n\setminus{W_n^\star}|=\mathcal O(n^{d-1})$ and $
|W_n^\star\setminus{W_n^\star}_{\ominus
2\varepsilon}|=\mathcal
O(n^{d-1+\alpha})$, since $\alpha<1$. Hence, 
%
\begin{equation}\label{ekors1}
\E\bigl(\bigl(A_{n}-A_n^\star\bigr)\Par\bigr)=
\mathcal O \bigl(n^{d-1+\alpha}\bigr) +\mathcal O \bigl(n^{d-1}\bigr) =
\mathcal O \bigl(n^{d-1+\alpha}\bigr)
\end{equation}
and
%
\begin{equation}\label{ekors2}
\E\bigl(b_{n}-b_n^\star\bigr)= \mathcal O
\bigl(n^{d-1+\alpha}\bigr) +\mathcal O \bigl(n^{d-1}\bigr) = \mathcal O
\bigl(n^{d-1+\alpha}\bigr).
\end{equation}
Since $h_\varepsilon^{(n)}$ has support included in $W_n^\star$,
Corollary \ref{corest} gives $\E(A_n^\star\Par+
b_n^\star)=0$. Combining this with (\ref{ekors1})--(\ref{ekors2}) gives
the result of (II).

(III): From Lemmas \ref{lemdomain}--\ref{lemvar}, (ii), and (\ref
{eqhen}), we get
\[
\Var T_1 =\mathcal O \bigl(\bigl|W_n^\star
\setminus{W_n^\star}_{\ominus
2\varepsilon}\bigr| \bigr) = \mathcal O
\bigl(n^{d-1+\alpha}\bigr)
\]
and
\[
\Var T_2 = \mathcal O \bigl(\bigl|W_n\setminus{W_n^\star}\bigr|
\bigr) = \mathcal O\bigl(n^{d-1}\bigr),
\]
which leads to
%
\begin{equation}
\label{eqprob1} \Var\bigl( |W_n|^{-1/2}
\bigl(A_{n}-A_n^\star\bigr)\Par\bigr) = \mathcal
O \biggl( \frac{n^{d-1+\alpha}}{n^d} \biggr) = \mathcal O \bigl
({n^{\alpha
-1}}
\bigr).
\end{equation}
In the same way, we derive
\[
\Var T_1^\prime= \mathcal O \bigl(\bigl|W_n^\star
\setminus{W_n^\star}_{\ominus2\varepsilon}\bigr| \bigr)=\mathcal O
\bigl(n^{d-1+\alpha}\bigr)
\]
and
\[
\Var T_2^\prime=\mathcal O \bigl(\bigl|W_n\setminus
W_n^\star\bigr| \bigr)=\mathcal O\bigl(n^{d-1}\bigr),
\]
which leads to
%
\begin{equation}
\label{eqprob2} \Var\bigl( |W_n|^{-1/2}
\bigl(b_{n}-b_n^\star\bigr) \bigr) = \mathcal{O}
\bigl( {n^{\alpha-1}} \bigr).
\end{equation}
Combining (\ref{eqprob1})--(\ref{eqprob2}) with Chebyshev's
inequality completes the proof of (III).
\end{pf}

Finally, we turn to the proof of (a)--(c) in Theorem \ref{thmconv}. %

(a): With probability one, by (I) in Lemma \ref{lemZn},
$|W_n|^{-1}(A_{n}-S_n) \geq-|W_n|^{-1} S_n/2$ for all sufficiently
large $n$, and so by (iii),
%
\begin{equation}
\label{eqlbAn} \frac{A_{n}}{|W_n|} \geq\frac{S_n}{2|W_n|} \geq\frac{I_0}2
\end{equation}
for all sufficiently
large $n$. Thereby, (a) 
is obtained.

(b): With probability one, for $n$ large enough,
we can write $|W_n|^{-1}A_{n}(\ParE-\Par)= -
|W_n|^{-1}\*(A_{n}\Par+b_{n})$, and
by (\ref{eqlbAn}), $\|(|W_n|^{-1}A_n)^{-1}\| \leq2/\mu_{\min}$ where
$\mu_{\min}$ is the smallest eigenvalue of~$I_0$. Combining this with
(a) in Theorem \ref{thmconv}, with probability one, for $n$
large enough, we obtain
\begin{eqnarray*}
\| \ParE- \Par\| &=& \bigl\| \bigl(|W_n|^{-1}{A_n}
\bigr)^{-1} |W_n|^{-1}(A_{n}
\Par+b_{n}) \bigr\|
\\
&\leq& \frac{2}{\mu_{\min}} \bigl\| |W_n|^{-1}(A_{n}
\Par+b_{n}) \bigr\| .
\end{eqnarray*}
The right-hand side of this inequality converges almost surely to
zero, cf. Lemma \ref{lemZn}. Thereby (b) 
follows.

(c): For a function $\psi\dvt\RR^d\to\RR$ and a bounded Borel set
$\Delta\subset\RR^d$, define 
%
\begin{equation}
\label{eqVDelta} V_\Delta(\psi) = \int_\Delta\psi(u)
\psi(u)^\top\rho(u)\dd u + \int_\Delta\int
_\Delta\psi(u_1)\psi(u_2)^\top
Q_2(u_1,u_2)\dd u_1\dd
u_2
\end{equation}
provided the integrals exist (are finite).
Observe that $\Sigma_n = V_{W_n}(f^{(n)}_\Par)$ and
$\Sigma_n^\star=V_{W_n^\star}(f^{(n)}_{\Par,\varepsilon})$ where
\[
f^{(n)}_{\Par,\varepsilon}(u)=h^{(n)}_{\varepsilon}(u) \div
z(u)^\top\Par+\div h_{\varepsilon}^{(n)}(u).
\]
We decompose the proof of (c) 
into three steps.

\textit{Step} 1. Assuming $\Sigma_n^\star\geq I_0>0$ for some positive
definite matrix $I_0$ and for all $n$ large enough, we prove that
%
\begin{equation}
\label{e12} {\Sigma_n^\star}^{-1/2}
\bigl(A_n^\star\Par+b_n^\star\bigr)
\xrightarrow{d} \mathcal{N}(0,I_p) \qquad\mbox{as $n\rightarrow\infty$}.
\end{equation}
We have
\[
A_n^\star\Par+ b_n^\star= \sum
_{i \in\mathcal I_n} Y_i^{(n)} \qquad\mbox{with }
Y_i^{(n)}=\sum_{u \in X_{C_i}}
f^{(n)}_{\Par,\varepsilon}(u).
\]
For any $n\geq1$ and any $i \in\mathcal{I}_n$, $Y_i^{(n)}$ has zero
mean, and by (iv),
\[
\sup_{n\geq1} \sup_{i\in\mathcal I_n}\E\bigl(
\bigl\| Y_i^{(n)} \bigr\| ^{2+\delta} \bigr)=\mathcal O(1).
\]
This combined with (v) and the assumption on $\Sigma_n^\star$, allows
us to invoke Kar\'aczony (\cite{karaczony06}, Theorem~4), which is a
central limit
theorem for a triangular array of random fields, which in turn is
based on Guyon (\cite{guyon91}, Theorem 3.3.1). Thereby (\ref{e12})
is obtained.

\textit{Step} 2. We prove that
%
\begin{equation}
\label{e13} |W_n|^{-1}\bigl(\Sigma_n -
\Sigma_n^\star\bigr) \to0 \qquad\mbox{as $n\rightarrow\infty$}.
\end{equation}
Using the notation (\ref{eqVDelta}), we have
%
\begin{equation}
\label{eqdiffS} \Sigma_n-\Sigma_n^\star=
V_{{W_n^\star}_{\ominus2\varepsilon
}}\bigl(\zeta^{(n)}\bigr) + V_{W_n\setminus{W_n^\star}_{\ominus
2\varepsilon}}\bigl(
\zeta^{(n)}\bigr),
\end{equation}
where 
%
\begin{equation}
\zeta^{(n)}(u_1,u_2) =
f^{(n)}_\Par(u_1)f^{(n)}_\Par(u_2)^\top-
f_{\Par,\varepsilon}^{(n)}(u_1)f_{\Par,\varepsilon}^{(n)}(u_2)^\top
,\qquad u_1,u_2\in\RR^d.
\end{equation}
By (ii) and (\ref{eqhen}), every entry of $\zeta^{(n)}(u_1,u_2)$
vanishes if $u_1,u_2\in{W_n^\star}_{\ominus2\varepsilon}$, and its
numeric value is
bounded by a constant if $u_1,u_2\in W_n$. Therefore, we can apply
similar arguments as used in the proof of Lemma \ref{lemvar} to
conclude that
\[
|W_n|^{-1} \bigl|\bigl(\Sigma_n-
\Sigma_n^\star\bigr)_{jk}\bigr| = |W_n|^{-1}
\bigl(V_{W_n\setminus{W_n^\star}_{\ominus2\varepsilon}}\bigl(\zeta
^{(n)}\bigr) \bigr)_{jk} =
\mathcal{O} \biggl( \frac{|W_n^\star\setminus{W_n^\star
}_{\ominus
2\varepsilon}|}{|W_n|} \biggr) = \mathcal{O} \bigl(n^{\alpha-1}
\bigr),
\]
which leads to the verification of (\ref{e13}).


\textit{Step} 3.
From (vi) and (\ref{e13}), we see that with probability one,
$\Sigma_n^\star$ is invertible for all sufficiently large $n$,
which allows us to write
%
\begin{eqnarray}
\Sigma_n^{-1/2} S_n (\ParE-\Par) &=& -
\Sigma_n^{-1/2}(A_{n}\Par+b_{n})
\nonumber
\\
\label{eqclt1}
&=&- \Sigma_n^{-1/2}\bigl(\bigl(A_{n}-A_{n}^\star
\bigr)\Par+b_{n}-b_{n}^\star\bigr)
\\
\label{eqclt2}
&&{} + \bigl( \Sigma_n^{-1/2}-\bigl(\Sigma_n^\star
\bigr)^{-1/2} \bigr) \bigl(A_{n}^\star
\Par+b_{n}^\star\bigr)
\\
&&{}+ \bigl(\Sigma_n^\star\bigr)^{-1/2}
\bigl(A_{n}^\star\Par+b_{n}^\star\bigr).
\nonumber
\end{eqnarray}
From (\ref{e12}) and Slutsky's lemma, we obtain that
(\ref{eqconvDist}) will be true if we manage to prove that the two
terms (\ref{eqclt1}) and (\ref{eqclt2}) converge towards zero in
probability as $n\to\infty$. Let $U_1$ and $U_2$ denote these two
terms. Let $M_n=\Sigma_n^\star/|W_n|$.
For $n$ large enough, we have $\|M_n^{-1}\|\leq2/{\lambda_{\min}}$,
so $\|M_n^{-1/2}\|\leq2/\sqrt{\lambda_{\min} }$, where
$\lambda_{\min}$ is the smallest eigenvalue of $I_0'$ in (vi), and
there exists a constant $C$ such that $\max(\|M_n^{1/2}\|,\|M_n\|)\leq
C$. On the first hand, we note that
\[
\| U_1 \| \leq\frac{2}{\sqrt{\lambda_{\min}}} \bigl\| |W_n|^{1/2}
\bigl( \bigl(A_{n}-A_{n}^\star\bigr)\Par
+b_{n}-b_{n}^\star\bigr) \bigr\|,
\]
which from (III) in Lemma \ref{lemZn} leads to $U_1\xrightarrow{P}
0$ as $n\to\infty$. On the other hand, we have
%
\begin{equation}
\label{eqU2} U_2 = \bigl( \Sigma_n^{-1/2}\bigl(
\Sigma_n^\star\bigr)^{1/2}-I_p \bigr)
\bigl(\Sigma_n^\star\bigr)^{-1/2}
\bigl(A_{n}^\star\Par+ b_{n}^\star\bigr).
\end{equation}
Since $\|(\Sigma_n / |W_n|)^{-1}\|$ is bounded, we derive from (\ref
{e13}) that
\[
\biggl( \frac{ \Sigma_n}{|W_n|} \biggr)^{-1} \biggl(
\frac{\Sigma
_n-\Sigma
_n^\star}{|W_n|} \biggr) = I_p-{\Sigma}_n^{-1}
\Sigma_n^\star\rightarrow0,
\]
which also leads to ${\Sigma}_n^{-1/2}
(\Sigma_n^\star)^{1/2}\to I_p$. Combining (\ref{e12}) and
(\ref{eqU2}) with Slutsky's lemma, convergence in
probability to zero of $U_2$ is deduced. The proof of
Theorem \ref{thmconv} is thereby completed.

\section{The \textsc{vare} for model 2} \label{appmodel2}

For specificity and simplicity, consider
the setting of Section \ref{secsimd2} when $h(u)=\div
z(u)$ and model~2 is assumed. Then a
straightforward calculation leads to the following simple expression
for the \textsc{vare}:
\[
\ParE= \pmatrix{ \displaystyle \sum\cos^2(4\uppi  u_1) &
\displaystyle \sum
\cos(4\uppi  u_1)\cos(4\uppi  u_2)
\vspace*{2pt}\cr
\displaystyle \sum
\cos(4\uppi  u_1)\cos(4\uppi  u_2) & \displaystyle \sum
\cos^2(4\uppi  u_2)}^{-1} \pmatrix{ \displaystyle \sum
\sin(4\uppi  u_1)
\vspace*{2pt}\cr
\displaystyle \sum\sin(4\uppi  u_1)},
\]
where $u=(u_1,u_2)\in\RR^2$ and $\sum=\sum_{u\in \mathbf{X}_{W_n}}$. In the
sequel, we discuss the conditions (i)--(vi) specified in
Section \ref{secconditions}.

Conditions (i), (iv), and (v) are discussed in
Section \ref{secdisccond} and are satisfied under the setting of
Section \ref{secsimd2}. Condition (ii) is obviously satisfied for
model 2. Below we focus on condition (iii) as condition
(vi) can be checked using similar ideas.

According to the discussion in Section \ref{secdisccond}, we only need
to verify that $|W_n|^{-1}\widetilde S_n \geq I_0$ where
\[
\widetilde S_n = \int_{W_n} \div z(u) \div
z(u)^\top\dd u.
\]
Let $C_i$ denote the unit cube centered at $i \in\mathcal{I}_n$ where
\[
\mathcal I_n =\bigl\{ (j,k)\dvt  j,k\in\{ -n/2,\ldots,-1/2,1/2,\ldots,n/2\}\bigr\}.
\]
Then $W_n=[-n,n]^2=\bigcup_{i\in\mathcal I_n} C_i$. Let $\eta>0$. There
exists a non-negative real-valued continuous function $f$ such that
$f(\eta)\to0$ as $\eta\to0$, and such that for any $i=(i_1,i_2) \in
\mathcal I_n$ and any $u=(u_1,u_2) \in b((i_1,i_1-3/8),\eta)$
\[
\bigl|\cos(4\uppi  u_1) -1 \bigr| \leq f(\eta) \quad\mbox{and}\quad \bigl|\cos(4\uppi
u_2)\bigr|\leq f(\eta).
\]
Therefore, for any $u \in b((i_1,i_1-3/8),\eta)$ and $y\in\RR^2
\setminus\{0\}$, whenever $\eta$ is sufficiently small,
\begin{eqnarray*}
&&
y^\top\div z(u) \div z(u)^\top y
\\
&&\quad= 16\uppi ^2 \bigl( y_1^2
\cos^2(4\uppi  u_1)+ 2y_1y_2 \cos(4
\uppi  u_1)\cos(4\uppi  u_2) + y_2^2
\cos^2(4\uppi  u_2) \bigr)
\\
&&\quad\geq16\uppi ^2 \bigl( y_1^2 \bigl(1-f(\eta)
\bigr) -2 |y_1 y_2| f(\eta)^2
-y_2^2 f(\eta) \bigr) \geq8 \uppi ^2
y_1^2.
\end{eqnarray*}
Thus, for sufficiently small $\eta$,
\[
y^\top\widetilde S_n y = \sum
_{i \in\mathcal I_n} \int_{C_i} y^\top\div
z(u) \div z(u)^\top y \dd u\geq8\uppi ^2y_1^2
\bigl(\uppi \eta^2\bigr) |\mathcal I_n| = c
|W_n|
\]
with $c=8\uppi ^3 y_1^2 \eta^2>0$. This implies that
$|W_n|^{-1}\widetilde
S_n \geq c J_2$ where $J_2$ is the $2\times2$ identity matrix.
\end{appendix}

\section*{Acknowledgments}

This research was initiated when J.-F. Coeurjolly was a Visiting
Professor at Department of Mathematical Sciences, Aalborg University,
February--July 2012, and he thanks the members of the department for
their kind hospitality. The research of J.-F. Coeurjolly was also
supported by Joseph Fourier University of Grenoble (project
``SpaComp''). The research of J. M\o ller was supported by the Danish
Council for Independent Research|Natural Sciences, Grants
09-072331 (``Point process modelling and statistical inference'') and
12-124675 (``Mathematical and statistical analysis of spatial data''),
and by the Centre for Stochastic Geometry and Advanced Bioimaging,
funded by a grant from the Villum Foundation. Both authors were
supported by l'Institut Fran\c{c}ais du Danemark.



\printhistory

\end{document}